\documentclass[11pt]{article}
\usepackage{amsmath,amsfonts,amssymb,color}
\usepackage{graphicx}
\usepackage[mathscr]{euscript}
 \usepackage{setspace}
\newcommand{\ep}{\varepsilon}
\newcommand{\al}{\alpha}
\newcommand{\ov}{\overline}

\newcommand{\R}{\mathbb R}

\newcommand{\eq}{equilibrium }

\begin{document}   
\title{ Diffusive Heat Transport in Budyko's Energy Balance Climate Model with a Dynamic Ice Line}
\author{James Walsh\\Department of Mathematics, Oberlin College}
\date{July 2016}

\maketitle

\vspace{.5in}
\noindent{\bf Abstract.} \ M. Budyko and W. Sellers independently introduced seminal  energy balance climate models in 1969,  each with a goal of investigating the role played by positive ice albedo feedback in climate dynamics.  In this paper we replace the relaxation to the mean  horizontal  heat transport mechanism used in the models of Budyko and Sellers with diffusive heat transport. We couple the resulting surface temperature equation with an equation for movement of the edge of the ice sheet (called the ice line), recently introduced by E. Widiasih. We apply the spectral method to the temperature-ice line system   and   consider   finite approximations. We prove there exists a stable \eq solution with a small ice cap, and an unstable \eq solution with a large ice cap, for a range of parameter values. If the diffusive transport is too efficient, however, the small ice cap disappears and  an ice free Earth becomes a limiting state. In addition, we analyze a variant of the coupled diffusion equations appropriate as a model for extensive glacial episodes in the Neoproterozoic Era. Although the model equations are no longer smooth due to the existence of a switching boundary, we prove there exists a unique stable equilibrium solution with the ice line in tropical latitudes, a climate event known as a  Jormungand or Waterbelt state. As the systems introduced here contain variables with differing time scales, the main tool used in the analysis is geometric singular perturbation theory.

\bigskip



\section{Introduction}

Techniques from dynamical systems theory are playing an ever increasing role in the study  of low-order climate models. From finite-dimensional approximations of the physics-based model PDE  (\cite{Lor84}, \cite{maas},  \cite{shil}, \cite{swart}, \cite{veen2003}), to models investigating low-frequency variability (\cite{broer2008}, \cite{broer2011}, \cite{roeb}, \cite{vv2003}) and mixed-mode oscillations \cite{rob}, to new glacial cycle models (\cite{ash}, \cite{wwhm}), and to the direct analysis of reaction diffusion PDE defined on the 2-sphere (\cite{diaz2006},  \cite{diaz2002},   \cite{diaz2008}, \cite{hetzer}, Part II of \cite{diazbook}),  mathematical analysis plays an important role in climate modeling. Conceptual climate models in particular lend themselves to a mathematical treatment,  with their use continuing to contribute to a greater  understanding of planetary climate \cite{held}.

Energy balance models (EBM) of climate are a type of conceptual model that focuses on major climate components and their interactions. Incoming solar radiation (or  {\em insolation}), a planet's outgoing longwave radiation, heat transported by the atmosphere and oceans, and planetary albedo (or reflectivity) provide examples of such components. Early and important models of M. Budyko \cite{bud} and W. Sellers \cite{sell} were introduced to investigate the effect of ice albedo on temperature. An expanding ice sheet increases albedo, for example, decreasing the absorption of insolation and thereby decreasing temperature, so that the ice sheet expands further. Conversely, a retreating ice sheet lowers planetary albedo, increasing temperature and further reducing the size of the ice sheet. Might this ``positive ice albedo feedback" have been the cause of snowball Earth events in our planet's past, during  which it is believed glaciers advanced to the equator? Would a retreating ice sheet inevitably result in an ice free planet?

The models of Budyko and Sellers consider  average annual surface temperature in latitudinal zones, so that there is one spatial dimension. The horizontal (i.e., meridional) transport of heat was modeled in \cite{bud} via the relaxation of zonal temperature to the global mean temperature (equation \eqref{relax}, section 2.1). Budyko's analysis indicated that a bifurcation, or climate tipping point, occurs: If the annual global average insolation decreases by roughly 5\%, a small stable ice cap is lost and the climate  system heads to a snowball Earth.

It should be noted the analyses of early EBM, such as Budyko's equation, focused on equilibrium solutions and their sensitivity to changes in parameters. There is no mechanism, for example, by which the edge of the ice sheet (called the {\em ice line}) is allowed to move in response to changes in temperature. In \cite{wid}, E. Widiasih coupled zonal surface temperature with a dynamic ice line (with heat transport as in equation \eqref{relax}). From a dynamical systems perspective, this coupled model was infinite-dimensional, comprised of a  function space (temperature) crossed with an  interval (the ice line). Using Hadamard's graph transform method, Widiasih proved the existence of an invariant, locally exponentially attracting, one-dimensional manifold on which the coupled system was well-approximated by a single ODE. In particular, she proved the existence of a stable equilibrium solution with a small ice cap, and an unstable equilibrium solution with a large ice cap. Widiasih's work further indicated  the presence of a sufficiently large ice sheet would trigger a runaway snowball Earth event, while a very small ice cap would expand and limit on the stable small ice cap position. In particular, an ice sheet retreating to an ice free Earth state is  not a possibility in this coupled model, assuming meridional heat transport is given by  equation \eqref{relax}.

There is ample evidence indicating that ice sheets flowed into the ocean in tropical latitudes during two extensive glacial periods in the Neoproterozoic Era (roughly 630 and 715 million years ago; see \cite{ray} and references therein). This has led many to posit the existence of a snowball Earth during these times. Evidence also exists, however, indicating that photosynthetic organisms and more biologically complex species of sponge thrived before and after these glacial episodes \cite{abb}. This has led others to believe that a stable climate state with ice sheets in tropical latitudes, yet with a strip of open ocean water about the equator, existed. Abbot et al called  this  a {\em Jormungand} climate state in \cite{abb}, proposing an adjustment in the latitude-dependent albedo function as a mechanism  by which this stable state might exist  (and finding  a Jormungand state in an idealized, appropriately parametrized general circulation model). A variant of this albedo function was incorporated into the coupled Budyko--ice line model in \cite{walshwid}, in which it was shown the infinite-dimensional dynamical system has a stable equilibrium solution with the ice line within 20$^\circ$ of the equator, again assuming horizontal heat transport is modeled as in equation  \eqref{relax}.

Returning to Budyko's  original albedo function  and continuing to use the relaxation to the mean heat transport term \eqref{relax}, a much simpler quadratic approximation to the   coupled temperature--ice line model was introduced in \cite{dickesther}. This approach resulted in a reduction to a two-dimensional system of ODE that was analyzed using more elementary dynamical systems techniques. As with the original Budyko--ice line model, a small stable ice cap and a larger unstable ice cap were shown to exist. A similar quadratic approximation approach  in the case of the Jormungand temperature--ice line model was followed in \cite{walshrack}, where it was shown the quadratic approximation yielded an equilibrium solution with the ice edge in tropical latitudes. Notably, and due to the relative complexity of the Jormungand albedo function, the analysis in \cite{walshrack} required tools from the theory of discontinuous differential equations, as well as geometric singular perturbation theory. The latter theory   will play an important role in the present work as well.

In this paper the relaxation to the mean heat transport term is replaced by diffusive heat transport  in Budyko's surface temperature equation \eqref{bud}, leading to a variant of Budyko's   model discussed in several previous works, including \cite{heldsuarez}, \cite{linnorth}, \cite{north2}, \cite{north},  \cite{northetal}, and Chapter 10 in \cite{ghilbook}. In each of these previous studies  the ice line remains fixed in time. Budyko's zonal temperature PDE with diffusive heat transport is coupled here with Widiasih's dynamic ice line equation, and the resulting system is analyzed using the spectral method, most significantly in the case of the extensive glacial episodes of the Neoproterozoic Era. This approach produces a continuous vector field that is, however,  nonsmooth on a  ``switching boundary" \cite{dibernardo}. We show the vector field is locally Lipschitz, guaranteeing uniqueness of solutions. Geometric singular perturbation theory is used to show the existence of an attracting Jormungand state, with ice line in the tropics, in the case of diffusive heat transport. 

For completeness, we apply the spectral method to the diffusive heat transport variant of the standard albedo Budyko--ice line coupled model as well. We also revisit  work in \cite{dickesther} and \cite{walshrack} that makes use of the relaxation to the mean heat transport, albeit incorporating higher-order approximations of the temperature function to derive the governing equations. The results in this paper might then be viewed as completing the analysis of the Budyko--ice line coupled model begun in \cite{wid}, \cite{walshwid}, \cite{dickesther} and \cite{walshrack}, an analysis informed by a dynamical systems perspective.

In the following section we present the Budyko--ice line coupled model. In section 3 we apply the spectral method to the diffusive heat transport model when using Budyko's standard albedo function.  Section 4 contains an analysis of the Jormungand temperature--ice line model in the case of diffusive heat transport. In  section 4 we prove   uniqueness of solutions for the resulting nonsmooth system of ODE, as well as the existence of a stable equilibrium with ice cap extending to the tropics. The Budyko and Jormungand ``relaxation to the mean" horizontal heat transport models are analyzed in section 5, albeit with higher order approximations than those used in \cite{dickesther} and \cite{walshrack}. We summarize this work in the concluding section.

\section{The coupled temperature--ice line model}

\subsection{Budyko's conceptual model}

The EBM introduced by Budyko focuses on the annual mean surface temperature in latitudinal zones. As mentioned above, analyses of EBM historically centered on equilibrium solutions and their sensitivity to changes in parameters (typically, annual globally averaged insolation). Energy balance was achieved when
\begin{equation}\label{balance}
E_{in}-E_{out}=E_{tran}
\end{equation}
at a given latitude, where $E_{in}$ represents  insolation (modulo the albedo), $E_{out}$ is the outgoing longwave radiation (OLR), and $E_{tran}$ is a meridional heat transport term. 

One assumes a symmetry about the equator in Budyko's model, so that latitude extends  from $\theta=0^\circ$ to  $\theta=90^\circ$. Additionally, one introduces the variable $y=\sin\theta$ (the area of an infinitesimal band at ``latitude" $y$ has area proportional to $dy$),  so that  $y=0 $ at the equator and $y=1$ at the North Pole. Any function dependent upon $y$ is thus an even function of $y$, given the symmetry assumption.

The $E_{in}$-term has the form 
\begin{equation}\label{Ein}
E_{in}=Qs(y)(1-\al(y,\eta)),
\end{equation}
where $Q$ is the annual global mean insolation in Watts per meter squared (W/m$^2$), and $s(y)$ provides for the distribution of insolation as a function of latitude (see Figure 1). Thus, $Qs(y)$ is the average annual insolation at latitude $y$.

The explicit formulation
\begin{equation}\label{realsofy}
s(y)=\frac{2}{\pi^2}\int^{2\pi}_0 \sqrt{1-(\sqrt{1-y^2} \, \sin\beta\cos\gamma-y\cos\beta)^2} \, d\gamma,
\end{equation}
where $\beta$ denotes the obliquity (or tilt) of the Earth's spin axis, was derived in \cite{dickclar}. We set  $\beta=23.5^\circ$, the current value of Earth's obliquity. Anticipating the use of the spectral method to follow, one can express \eqref{realsofy} in terms of the even Legendre polynomials $p_{2n}(y)$, writing
\begin{equation}\label{fullsLeg}
s(y)=\sum^\infty_{n=0} s_{2n}p_{2n}(y),
\end{equation} 
where $s_{2n}$ is given by
\begin{equation}\label{s2ns}
s_{2n}=(4n+1)\int^1_0 s(y)p_{2n}(y)dy.
\end{equation}
One finds the insolation components fall off rapidly, to the extent the quadratic
\begin{equation}\label{squad}
s(y)=s_0p_0(y)+s_2p_2(y), \ s_0=1,  \ s_2=-0.477,
\end{equation}
$ p_0(y)=1 \mbox{ and } p_2(y)=\frac{1}{2}(3y^2-1)$, approximates equation \eqref{realsofy} uniformly to within 3\% error (see Figure 1). 

Including higher-order terms from \eqref{fullsLeg} in the analysis did not alter the nature of the results in all of our investigations. Hence we replace equation \eqref{realsofy} with equation \eqref{squad} in all that follows, effectively setting $s_{2n}=0$ for $n>1$ in \eqref{fullsLeg} (an  exception to this approach will occur  in section 5, where at times we include higher-order $s_{2n}$-terms).

As the function $s(y)p_{2n}(y)$ occurs frequently below,   for ease of notation we set
 $$q_{2n}(y)=s(y)p_{2n}(y)$$ in all that follows.

The function $\al(y,\eta)$ in equation \eqref{Ein} represents the albedo, which we assume depends  upon the position of the ice line $y=\eta$. The  albedo function used by Budyko  has the form
\begin{equation}\label{budalb}
\al(y,\eta)=
 \begin{cases}
\al_1, &\text{if  \ } y<\eta \\
\al_2, &\text{if  \ } y>\eta,\\
\end{cases} \quad \quad
\end{equation}
with $\al_1<\al_2 \, $  \cite{bud}.  The larger value $\al_2$ represents the albedo of ice, with ice assumed to exist at all latitudes poleward of the ice line. The value $\al_1$ is the correspondingly  lower albedo of the  ice free surface equatorward of $\eta$. Each of $s(y)$ and $\al(y,\eta)$ is dimensionless.

\begin{figure}[t!]
\begin{center}
\includegraphics[width=5in,trim = 1.2in 7.1in 1.2in  1in, clip]{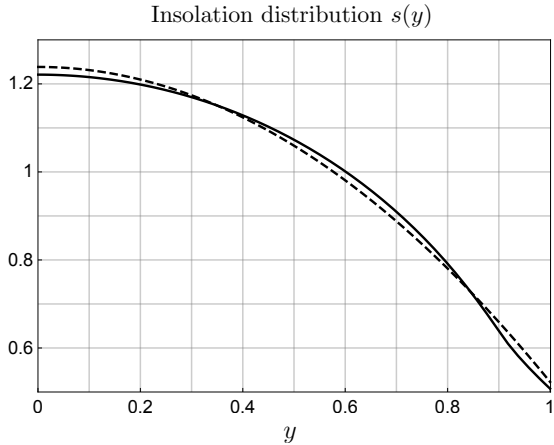}\\
\parbox{5in}{\caption{
{\em Solid}: Plot of equation \eqref{realsofy} with $\beta=24.5^\circ$. 
 {\em Dashed}: Quadratic approximating  equation \eqref{squad}. }}
 \end{center}
\end{figure}

The Earth radiates energy into space at longer wavelengths; following \cite{bud}, we set
\begin{equation}\label{Eout}
E_{out}=A+BT(y,t),
\end{equation}
where $T(y,t)$  ($^\circ$C) is the annual mean surface temperature at latitude $y$, and $A$ (W/m$^2$) and $B$ (W/(m$^2$ $ \!  ^\circ$C) are determined empirically (see, for example, \cite{graves}). We note Sellers chose a formulation of the $E_{out}$-term comprised of  an empirically determined  scaling of the Stephan-Boltzmann Law for blackbody radiation \cite{sell}. 

The $E_{tran}$-term in equation \eqref{balance} models the heat flux resulting from horizontal redistribution by the circulation of the ocean and the atmosphere.
In each of \cite{bud} and \cite{sell}, the horizontal heat transport term was given by  a relaxation to the mean expression of the form
\begin{equation}\label{relax}
E_{tran}=C(\overline{T}-T),
\end{equation}
where $C>0$ is an empirical constant and
\begin{equation}\label{Tbar}
\overline{T}=\ov{T}(t)=\int^1_0 T(y,t) dy
\end{equation}
is the global mean surface temperature. A latitude at which the temperature is below the global average will gain heat flux from neighboring latitudes, while heat will be transferred away horizontally if the local temperature exceeds the global mean.

In other variants of the Budyko-Sellers EBM (\cite{linnorth}, \cite{north2}, \cite{north}, \cite{northetal}, for example), meridional heat transport is modeled by
\begin{equation}\label{diffuse}
E_{tran}=D\nabla^2 T=D\frac{\partial}{\partial y}(1-y^2)\frac{\partial T}{\partial y},
\end{equation}
where $D>0$ is an empirical constant; the form of the spherical diffusion operator \eqref{diffuse} arises when assuming  no radial or longitudinal dependence. This diffusive heat transport is akin to an eddy diffusion approach to dispersion by macroturbulence in the entire system \cite{linnorth}.

In using either equation \eqref{relax} or equation \eqref{diffuse}, note the net global effect of horizontal heat transport is zero as in each case the integral of $E_{tran}$ over the unit interval is zero.

The diffusive heat transport approach \eqref{diffuse} has the advantage that the Legendre polynomials $p_n(y)$ are   eigenfunctions of the spherical diffusion operator:
\begin{equation}\notag
\frac{d}{dy}(1-y^2)\frac{d}{dy}p_n(y)=-n(n+1)p_n(y).
\end{equation}
Hence the  use of Legendre polynomials provides for a powerful technique in the analysis to follow.

We  consider the PDE version of Budyko's model
\begin{equation}\label{bud}
R\frac{\partial T(y,t)}{\partial t}=Q s(y) (1-\al(y,\eta))-(A+BT)+E_{tran},
\end{equation}
where $E_{tran}$ is of the form \eqref{relax} or \eqref{diffuse}. The heat capacity of the Earth's surface $R$ \, (J/(m$^{2} \,  ^\circ$C)) is assumed to equal the global average value; the units on each side of equation \eqref{bud} are W/m$^2$.

\subsection{Widiasih's dynamic ice line}

Note the lack of any mechanism in Budyko's model by which the ice sheet (represented in \eqref{bud} by the ice line $\eta$) is allowed to expand or retreat as the local surface temperature varies. In the case $E_{tran}=C(\ov{T}-T)$, a dynamic equation for the ice line was introduced in \cite{wid} and coupled to equation \eqref{bud}, resulting in the integro-differential  system
\begin{subequations}\label{budwid}
\begin{eqnarray}
R\frac{\partial T}{\partial t} & =   & Qs(y)(1-\al(y,\eta))-(A+BT)+C(\ov{T}-T)  \label{budwidA}  \\
\frac{d\eta}{dt} & =   &  \ep (T(\eta,t)-T_c).\label{budwidB}
\end{eqnarray}
\end{subequations}
The parameter $\ep>0$  governs the relaxation time of the ice sheet.   $T_c$ is a critical temperature, meaning that if $T(\eta,t)<T_c$ ice forms and  the ice line moves toward the equator, whereas the ice line retreats toward the pole if $T(\eta,t)>T_c$.  

System \eqref{budwid} has been analyzed in \cite{wid} and \cite{walshwid} when  assuming the albedo varies smoothly across the ice line (in this case the model is said to be of ``Sellers-type"). Using a smooth variant of albedo function \eqref{budalb} and assuming $\ep $ satisfies certain bounds, Widiasih   proved that a stable equilibrium temperature profile--ice line pair $(T^*(y),\eta^*)$ exists with $\eta^*$ at roughly 70$^\circ$N (a small ice cap). The analysis in \cite{wid} centered on the use of Hadamard graph transforms. 

With a smooth variant of albedo function \eqref{Jalb}, chosen  to model the very cold world of the great glacial epochs of the Neoproterozoic Era (see section 4), Walsh and Widiasih proved the existence of a stable equilibrium for system \eqref{budwid}
with the ice line in tropical latitudes in \cite{walshwid} (with constraints placed on $\ep$). The choice of albedo function was motivated by \cite{abb}, in which such a stable equilibrium was found using an idealized, appropriately parametrized general circulation model, and subsequently dubbed a ``Jormungand" state. That the temperature--ice line coupled model \eqref{budwid}  can be adjusted to find a stable climate appropriate for current times, as well as a stable climate appropriate for the deepest of ice ages, points to the power and utility of this relatively simple coupled EBM.

The present work considers discontinuous albedo functions---in which case equation \eqref{bud} is said to be of ``Budyko-type"---and diffusive heat transport. In the following two sections we present the analysis of the coupled temperature--ice line model   assuming diffusive heat transport \eqref{diffuse} and albedo functions \eqref{budalb} and \eqref{Jalb}, respectively. For completeness, in section 5 we then  revisit and expand upon  the quadratic approximation approaches to analyzing system \eqref{budwid} presented in \cite{dickesther} and \cite{walshrack}.

\section{Diffusive heat transport: standard albedo}

In this section we consider the system
\begin{subequations}\label{budwiddiff}
\begin{eqnarray}
R\frac{\partial T}{\partial t} & =   & Qs(y)(1-\al(y,\eta))-(A+BT)+D\frac{\partial }{\partial y}(1-y^2)\frac{\partial T}{\partial y}  \label{budwiddiffA}  \\
\frac{d\eta}{dt} & =   &  \ep (T(\eta,t)-T_c),\label{budwiddiffB}
\end{eqnarray}
\end{subequations}
with albedo function $\al(y,\eta)$ given by \eqref{budalb}. 

We note equation \eqref{budwiddiffA} comes with the boundary conditions that the gradient of the temperature profile $T(y,t)$ equals zero at the equator and at the North Pole.

\subsection{The spectral method}

 Recall functions appearing in system \eqref{budwiddiff} are even functions of $y$, due to the assumed symmetry about the equator. We approximate system \eqref{budwiddiff} by expanding each function of $y$ in terms of the first $N+1$ even Legendre polynomials, arriving at a system of $N+2$ ODE as follows. Consider the expression
 \begin{equation}\label{TLeg}
 T(y,t)=\sum^N_{n=0} T_{2n}(t)p_{2n}(y),
 \end{equation}
 and  first note
 \begin{equation}\notag
 \frac{d p_{2n}}{d\theta}=  \frac{d p_{2n}}{dy} \cos\theta =0
 \end{equation}
 at the equator and North Pole. Hence expression \eqref{TLeg} satisfies the prescribed boundary conditions. 
 
 Substitute \eqref{TLeg} and equation \eqref{squad} into equation \eqref{budwiddiffA} to arrive at
 \begin{align}\label{substitute}
 R\sum^N_{n=0}\dot{T}_{2n}p_{2n}(y)&=Q \sum^N_{n=0} (s_{2n}-\ov{\al}_{2n})p_{2n}(y) - A -B\sum^N_{n=0}T_{2n}p_{2n}(y)\\\notag
 &-D\sum^N_{n=0}2n(2n+1)T_{2n}p_{2n}(y),
 \end{align}
  where the dot represents differentiation with respect to time, $s_0$ and $s_2$ are given in \eqref{squad} (recall we set  $s_{2n}=0$ for $n>1$),
  and
 \begin{align}\label{al2nbars}
 \ov{\al}_{2n}(\eta)&=(4n+1)\int^1_0 \al(y,\eta)q_{2n}(y) dy\\\notag
 &=(4n+1)\left(\al_1\int^\eta_0 q_{2n}(y) dy+\al_2\int^1_\eta q_{2n}(y)dy\right).\\\notag
 &=(4n+1)\left(\al_2\int^1_0 q_{2n}(y) dy-(\al_2-\al_1)\int^\eta_0 q_{2n}(y)dy\right)\\\notag
 &=\al_2s_{2n}-(4n+1)(\al_2-\al_1)\int^\eta_0 q_{2n}(y) dy.\notag
 \end{align}
Recalling $q_{2n}(y)=s(y)p_{2n}(y)$ is a polynomial of degree $2n+2$, 
we note for future reference $\ov{\al}_{2n}$ is a polynomial in $\eta$ of degree $2n+3$.
Also note that substituting $y=\eta$ into expression \eqref{TLeg} gives the temperature at the ice line
\begin{equation}\label{Tofetadiff}
T(\eta,t)=\sum^N_{n=0}T_{2n}(t)p_{2n}(\eta).
\end{equation}

Equating coefficients of $p_{2n}$ in equation \eqref{substitute} and simplifying yields an approximation of the infinite-dimensional system \eqref{budwiddiff} given by the $N+2$ ODE 
\begin{align} \label{diffsyst}
 R\dot{T}_0&=Q(s_0-\ov{\al}_0(\eta))-A-BT_0, \\ \notag
R\dot{T}_{2n}&=Q(s_{2n}-\ov{\al}_{2n}(\eta))-(B+2n(2n+1)D)T_{2n}, \ n=1, ... , N,\\  \notag
\dot{\eta}&=\ep\left(\sum^N_{n=0}T_{2n}p_{2n}(\eta)-T_c \right). \notag
\end{align}

Note system \eqref{diffsyst} can be rewritten
\begin{align}\label{GSPdiff}
\dot{T}_0&=- \frac{B}{R}(T_0-f_0(\eta)),\\\notag
\\[-.15in]\notag\dot{T}_{2n}&=-\frac{(B+2n(2n+1)D)}{R}(T_{2n}-f_{2n}(\eta)), \ n=1, ... , N,\\\notag
\dot{\eta}&=\ep\left(\sum^N_{n=0}T_{2n}p_{2n}(\eta)-T_c \right),\notag
\end{align}
where 
\begin{align}\label{f2ns}
f_0(\eta)&=\frac{1}{B}(Q(s_0-\ov{\al}_0(\eta))-A),\\\notag
\\[-.15in]\notag f_{2n}(\eta)&=\frac{1}{(B+2n(2n+1)D)}(Q(s_{2n}-\ov{\al}_{2n}(\eta)), \ n=1, ... , N.\notag
\end{align}
We analyze system \eqref{GSPdiff} via geometric singular perturbation theory.

\subsection{Fenichel's Theorem}

We assume the ice line moves slowly relative to the evolution of the temperature. Thus $\eta$ will be  a slow variable, while $T_0, ... , T_{2N}$ are fast variables. 
If $\ep=0$ then $\eta$ remains fixed, in which case $T_{2n}\rightarrow f_{2n}(\eta)$ exponentially fast for $n=0, 1, ... , N$. Hence there is an attracting  manifold of rest points when $\ep=0$. We thus analyze system \eqref{GSPdiff} via geometric singular perturbation theory, for which we now recall Fenichel's Theorem (\cite{fen1}, \cite{fen2}, \cite{jones}). The case in which the attracting invariant manifold when $\ep=0$ is the graph of a function is sufficient for our purposes, so  this is the setting we present here.

Consider a system of the form
\begin{align}\label{fastslow}
\dot{x} & =  f(x,u,\ep)\\\notag
\dot{u} &= \ep g(x,u,\ep),\notag
\end{align}
where $x\in \R^n, u\in\R^m$ and $\ep$ is a real parameter. Assume $f$ and $g$ are each $C^\infty$ on $\R^{n+m}\times I$, where $I$ is an open interval containing 0. Assume the set of rest points when $\ep=0$,
\begin{equation}\notag
\Lambda^0 = \{(x,u) : f(x,u,0)=0\}=\{(x,u) : x=h^0(u)\},
\end{equation}
is the graph of a $C^\infty$ function $x=h^0(u)$. Assume further that $\Lambda^0$ is normally hyperbolic with regard to system \eqref{fastslow} when $\ep=0$.

\vspace{0.1in}
\noindent
{\bf Theorem 3.1} \cite{fen1}
{\em If $\ep$ is positive and sufficiently small, there exists a manifold $\Lambda^\ep$ within $O(\ep)$ of $\Lambda^0$ and diffeomorphic to $\Lambda^0$. In addition, $\Lambda^\ep$ is locally invariant  under the flow of system \eqref{fastslow} and $C^r$ (including in $\ep$) for any $r<\infty$.}

\vspace{0.1in}
\noindent
{\bf Remark 1.} \ 
(i) In this setting, $\Lambda^\ep$ is the graph of a function $x=h^\ep(u)$ that is  $C^r$, in both $u$ and $\ep$, for any $r<\infty$. \\
(ii) Given that $u$ parametrizes the manifold $\Lambda^\ep$, the equation describing the dynamics on $\Lambda^\ep$ is
\begin{equation}\label{regtime}
\dot{u}=\ep g(h^\ep(u),u,\ep).
\end{equation}
(iii) Changing the time scale ($\tau=\ep t$), equation \eqref{regtime} can be rewritten
\begin{equation}\notag
 \frac{du}{d\tau}= g(h^\ep(u),u,\ep)=g(h^0(u),u,0)+O(\ep).
\end{equation}
We then have
\begin{equation}\notag
\dot{u}=\ep g(h^\ep(u),u,\ep)=\ep(g(h^0(u),u,0)+O(\ep)).
\end{equation}
(iv) The function  $x=h^0(u)$ (and $x=h^\ep(u)$ as well) is defined on a compact subset $K$ of $\R^n$. In the applications of Theorem 3.1 to follow, $K$ can be any compact subset of $\R$.

\subsection{Standard albedo model dynamics}

Returning to system \eqref{GSPdiff}, we see that when $\ep=0$ there is a globally attracting curve of rest points $\Lambda^0=\{(x,\eta) : x=h^0(\eta)\}$, where 
\begin{equation}\label{h0}
x=h^0(\eta)=(f_0(\eta), f_{2}(\eta), ... , f_{2N}(\eta)),
\end{equation}
with $f_{2n}$ as  in \eqref{f2ns}, $n=0, 1, ... , N$. By Fenichel's Theorem, system \eqref{GSPdiff} has an attracting invariant manifold for sufficiently small $\ep$, on which the dynamics are well approximated by the single ODE
\begin{equation}\label{buddiffhofeta}
\dot{\eta}=\ep \, \left(\sum^N_{n=0} f_{2n}(\eta)p_{2n}(\eta)-T_c\right) = \ep h(\eta).
\end{equation}
Recalling $f_{2n}$ is a polynomial of degree $2n+3$, we see $h(\eta)$ is a polynomial of degree $4N+3.$ To each value $\eta=\eta^*$ with $h(\eta^*)=0$ and $h^\prime(\eta^*)<0$, there then corresponds a stable equilibrium solution of \eqref{GSPdiff}, near $(h^0(\eta^*),\eta^*)$, assuming the ice line moves sufficiently slowly.

We use the following parameter values  \cite{dickesther}, chosen to align with the modern climate:
\begin{equation}\label{budparams}
Q=343, A=202, B=1.9, \al_1=0.32, \al_2=0.62, T_c=-10.
\end{equation}

With the above parameters set, the existence of a small stable ice cap then depends upon the value of the diffusion constant $D$. A larger $D$-value corresponds to more efficient heat transport from the warmer tropical region to northern latitudes, eliminating the possibility of any small ice cap (see Figure 2(a)). Smaller $D$-values induce a steeper equator-to-pole temperature gradient, and the possibility of a small stable ice cap exists (Figure 2(b)), the size of which grows as $D$ decreases (compare Figures 2(b)  and 3(b)). 

\vspace{0.1in}
\noindent
{\bf Proposition 1.} \ 
{\em Assume   parameter values as in \eqref{budparams}, and let $N=1$ and $D=0.35$. For sufficiently small $\ep$, system \eqref{diffsyst} has two equilibrium solutions for which the ice line lies in the interval $[0,1]$. The equilibrium solution with the larger ice cap is unstable, while the equilibrium solution with the smaller ice cap is stable.}

\vspace{0.1in}
\noindent
{\em Proof}. This result follows from Fenichel's Theorem and the analysis of the degree 7 polynomial \eqref{buddiffhofeta}
\begin{equation}\notag
h(\eta)=f_0(\eta)p_0(\eta)+f_2(\eta)p_2(\eta)-T_c
\end{equation}
(see Figure 2(b)). One can  show $h(\eta)$ has two real roots $\eta_1<\eta_2$ in $[0,1]$, with $h^\prime(\eta_1)>0$ and $h^\prime(\eta_2)<0$. System \eqref{diffsyst} then has a stable equilibrium solution well approximated by $(h^0(\eta_2),\eta_2)$ and an unstable equilibrium solution near $(h^0(\eta_1),\eta_1)$, with $h^0$ as in equation \eqref{h0}.  \hfill $_\square$

\begin{figure}[t!]
\begin{center} \hspace*{.2in}
\includegraphics[width=6in,trim = 1.7in 6.54in .7in  1in, clip]{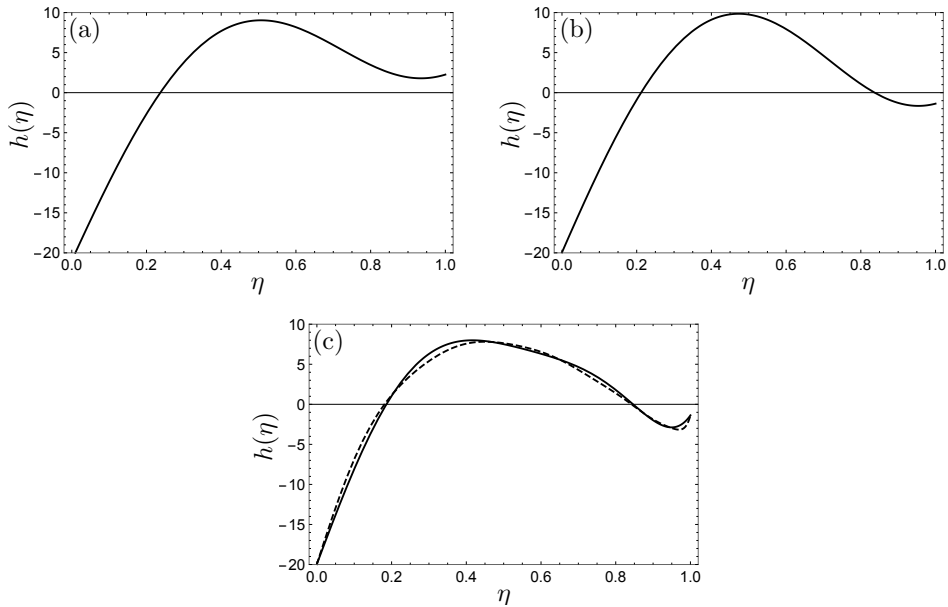}\\
\parbox{5in}{\caption{Plots of function \eqref{buddiffhofeta}. \ (a) $ N=1$ in \eqref{TLeg}, $ D=0.45$. \  (b) $N=1, D=0.35$. \ (c) D=0.35. {\em Solid}: $N=2$. {\em Dashed}: $N=5$.}}
 \end{center}
\end{figure}

\vspace{0.1in}
\noindent
{\bf Remark 2.} \ 
(i) Also plotted in Figure 2(c) are higher-order approximations ($N=2, 5$ in equation \eqref{TLeg}). As noted in \cite{northetal}, each succeeding term in \eqref{TLeg} contains system information pertaining to ever smaller spatial scales. Though not apparent in Figures 2(b) and 2(c), including higher-order terms serves to slightly reduce the size of the ice cap. We also see the quadratic approximation with $N=1$ suffices to capture the dynamics of system \eqref{diffsyst}.

\vspace{0.05in}
\noindent
(ii) It is interesting to note the lack of any unstable  \eq with   ice line  lying poleward of the stable ice line position when $D=0.35$ (Figure 2(b)). Previous analyses of (uncoupled) equation \eqref{budwiddiffA}  did find such a stable-unstable pair of equilibria, with the existence of the unstable ice cap known as {\em small ice cap instability}  (SICI) (see \cite{north84} for an excellent discussion of SICI and   Budyko's model with diffusive heat transport). The value $D=0.35$ is evidently sufficiently small in the coupled model \eqref{diffsyst} to preclude the existence of a small unstable ice cap. We note the SICI phenomenon does arise  when $D$ is increased to 0.394  (see Figure 3; compare with Figure 3 in \cite{north84}). Thus with more efficient meridional heat transport, the coexistence of a small stable ice cap and a ``stable" ice free Earth becomes a possibility, in contrast to the results for the relaxation to the mean heat transport model (\cite{wid}, \cite{dickesther};  also see Figure 6).

\begin{figure}[t!]
\begin{center} \hspace*{0.2in}
\includegraphics[width=6in,trim = 1.7in 6.5in .7in  1in, clip]{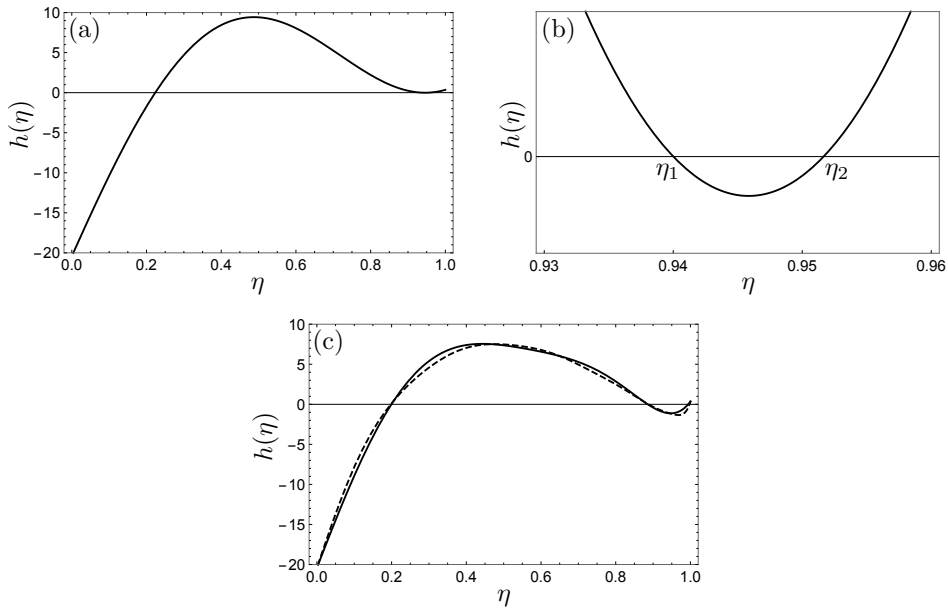}\\
\parbox{5in}{\caption{Plots of function \eqref{buddiffhofeta} with $D=0.394$. \ (a) and (b) $ N=1$ in \eqref{TLeg}. \  (c)  {\em Solid}: $N=2$ in \eqref{TLeg}. {\em Dashed}: $N=5$ in \eqref{TLeg}.}}
 \end{center}
\end{figure}

\vspace{0.05in}
\noindent
(iii) Integrating equation \eqref{TLeg} with respect to $y$ over the unit interval yields
\begin{equation}\notag
\ov{T}(t)=\int^1_0 T(y,t) dy=T_0(t),
\end{equation}
as $\int^1_0 p_{2n}(y) dy=0$ for $n>1$. Hence $T_0$ represents the global mean surface temperature, a quantity that varies with the ice line $\eta$ in the coupled model \eqref{diffsyst}. Given the evolution of $\eta$ depends on the $D$-value in system \eqref{diffsyst}, the global mean temperature depends upon $D$ as well. 
At an equilibrium point with $\eta=\eta^*$, the global mean temperature is $T^*_0=f_0(\eta^*)$ (recall equations \eqref{GSPdiff} and \eqref{f2ns}). With $N=1$ and $D=0.35$, the stable small ice cap is positioned near $\eta^*=0.837$ (57$^\circ$N), for which $T^*_0=10.9^\circ$C. In our present climate, $\eta^*\approx 0.940$ (70$^\circ$N) and $T^*_0\approx 15^\circ$C, so that the model ice cap is too large and the model world is too cold when $D=0.35$. Were one to require  the ice line for the small stable ice cap to lie
 near  $\eta^*=0.940$, one could compute   the corresponding $D$-value by solving $h(0.94)=f_0(0.94)+f_2(0.94)p_2(0.94)-T_c=0$ for $D$ (recall $f_2$ depends upon $D$). 
In this case one finds $D=0.394$, for which $T^*_0=14.6^\circ$C, a good approximation to our current climate.

\vspace{0.05in}
We now turn to diffusive heat transport and the spectral method when system \eqref{budwiddiff} is adjusted to model the great glacial episodes of the Neoproterozoic Era.

\section{Diffusive heat transport: Jormungand albedo}

\subsection{The Jormungand albedo function}

We incorporate an albedo function, modeled after a similar function introduced in \cite{abb} and appropriate for near snowball Earth events, into system \eqref{budwiddiff}. If the ice sheet were to extend down to tropical latitudes, what might serve to counteract positive ice albedo feedback  to the extent a runaway snowball Earth event never occurs, with  the ice line stabilizing in tropical latitudes?

A possible mechanism for this dynamic was presented in \cite{abb}. Using an idealized general circulation model (GCM) parametrized for the extensive ice ages of the Neoproterozoic Era, Abbot et al found that evaporation exceeded precipitation in a latitudinal band centered around 15-20$^\circ$. Any ice forming in these lower latitudes would have no (or significantly less) snow cover, and hence this ``bare" ice would have an albedo value $\al_i$ less than albedo value $\al_2$, but with $\al_i$  greater than the albedo value $\al_1$ of the surface south of the ice line.  This reduced albedo, relative to the albedo of snow covered ice, was sufficient to enable Abbot et al to find a stable \eq solution with the ice line in tropical latitudes when using the GCM in \cite{abb}.
This  stable ``Jormungand" state, with a narrow strip of ocean about the equator remaining ice free,   provides a possible explanation for the survival of photosynthetic organisms  through these tremendous ice ages.

Following \cite{abb}, we fix a latitude $y=\rho$ and assume any ice forming equatorward of $y=\rho$ is devoid of snow. In this scenario, a snowball Earth event would find bare ice from the equator to $y=\rho$, and snow covered ice poleward of $y=\rho$. If the ice line $\eta$ was positioned above $\rho$, the surface would have albedo $\al_1$ south of $\eta$ and albedo $\al_2$ at latitudes above $\eta$. This leads to the {\em Jormungand albedo function} \cite{walshrack}
\begin{equation} \label{Jalb}
\al(y,\eta)=
\begin{cases}
a(y,\eta), & \eta<\rho \\
b(y,\eta), & \rho\leq\eta,\\
\end{cases}
\end{equation}
where
\begin{equation} \label{andb}
\begin{array}{cc}
a(y,\eta)=
\begin{cases}
\al_1, & y<\eta \\
\al_i, & \eta<y<\rho\\
\al_2, & \rho<y,
\end{cases}  & \hspace*{.05in} \mbox{and}  \hspace*{.15in}
 b(y,\eta)=
\begin{cases}
\al_1, & y<\eta \\
\al_2, & y>\eta,
\end{cases}
\end{array}
\end{equation}
with  $\al_1<\al_i<\al_2$.

\subsection{The nonsmooth model equations}

The use of the function \eqref{Jalb} in system \eqref{budwiddiff} leads to a nonsmooth vector field as follows. Expanding $T(y,t)$ as in expression \eqref{TLeg} and proceeding as in section 3 again leads to the approximating system of $N+2$ ODE \eqref{diffsyst}. The functions $\ov{\al}_{2n}(\eta)$ appearing in system \eqref{diffsyst}, however, depend on the position of the ice line  relative to $y=\rho$. For $\eta\geq\rho$ we write $\ov{\al}_{2n}(\eta)=\ov{b}_{2n}$, and we have
\begin{align}\label{b2nbars}
\ov{\al}_{2n}(\eta)=\ov{b}_{2n}(\eta)&=(4n+1)\int^1_0b(y,\eta)q_{2n}(y)dy\\\notag
&=\al_2s_{2n}-(4n+1)(\al_2-\al_1)\int^\eta_0 q_{2n}(y) dy,\notag
\end{align}
precisely as in equation \eqref{al2nbars}. For $\eta<\rho$ we write $\ov{\al}_{2n}(\eta)=\ov{a}_{2n}$, and we compute
\begin{align}\label{a2nbars}
&\ov{\al}_{2n}(\eta)=\ov{a}_{2n}(\eta)=(4n+1)\int^1_0a(y,\eta)q_{2n}(y)dy\\\notag
 &=(4n+1)\left(\al_1\int^\eta_0 q_{2n}(y)dy+\al_i\int^\rho_\eta q_{2n}(y)dy+\al_2\int^1_\rho q_{2n}(y)dy\right)\\ \notag
 &=(4n+1)\left(\al_2\int^1_0 q_{2n}(y) dy-(\al_2-\al_i)\int^\rho_\eta q_{2n}(y)dy-(\al_2-\al_1)\int^\eta_0 q_{2n}(y)dy\right)\\ \notag
 &=\al_2s_{2n}-(4n+1)\left((\al_2-\al_i)\int^\rho_\eta q_{2n}(y)dy+(\al_2-\al_1)\int^\eta_0 q_{2n}(y)dy\right). \notag
\end{align}
We note $\ov{a}_{2n}(\eta)$ and $\ov{b}_{2n}(\eta)$ are each polynomials of degree $2n+3$.

Given albedo function \eqref{Jalb}, let  $V^-$ denote the vector field corresponding to model equations \eqref{diffsyst} when $\eta<\rho$, that is, $V^-$ is the vector field given by the system of $N+2$ ODE
\begin{align}\label{Vminus}
\dot{T}_0&=- \frac{B}{R}(T_0-f^-_0(\eta)),\\\notag
\\[-.5in]\notag\intertext{ \ \ $V^-$:}\notag
\\[-.7 in]\notag\dot{T}_{2n}&=-\frac{(B+2n(2n+1)D)}{R}(T_{2n}-f^-_{2n}(\eta)), \ n=1, ... , N,\\\notag
\dot{\eta}&=\ep\left(\sum^N_{n=0}T_{2n}p_{2n}(\eta)-T_c \right),\notag
\end{align}
where 
\begin{align}\label{fminus2ns}
f^-_0(\eta)&=\frac{1}{B}(Q(s_0-\ov{a}_0(\eta))-A),\\\notag
\\[-.15in]\notag f^-_{2n}(\eta)&=\frac{1}{(B+2n(2n+1)D)}(Q(s_{2n}-\ov{a}_{2n}(\eta)), \ n=1, ... , N.\notag
\end{align}
While $V^-$ is a smooth vector field on $\R^{N+2}$, in the model we  restrict $V^-$ to the set
\begin{equation}\notag
S^-=\{ (x,\eta) : x=(T_0, ... , T_{2N}), \ \eta<\rho\}.
\end{equation}

Similarly, let $V^+$ be the vector field given by the  model equations when $\rho<\eta$, that is, $V^+$ corresponds to the system of $N+2$ ODE
\begin{align}\label{Vplus}
\dot{T}_0&=- \frac{B}{R}(T_0-f^+_0(\eta)),\\\notag
\\[-.5in]\notag\intertext{ \ \ $V^+$:}\notag
\\[-.7in]\notag\dot{T}_{2n}&=-\frac{(B+2n(2n+1)D)}{R}(T_{2n}-f^+_{2n}(\eta)), \ n=1, ... , N,\\\notag
\dot{\eta}&=\ep\left(\sum^N_{n=0}T_{2n}p_{2n}(\eta)-T_c \right),\notag
\end{align}
where 
\begin{align}\label{fplus2ns}
f^+_0(\eta)&=\frac{1}{B}(Q(s_0-\ov{b}_0(\eta))-A),\\\notag
\\[-.15in]\notag f^+_{2n}(\eta)&=\frac{1}{(B+2n(2n+1)D)}(Q(s_{2n}-\ov{b}_{2n}(\eta)), \ n=1, ... , N.\notag
\end{align}
While $V^+$ is a smooth vector field on $\R^{N+2}$, we restrict  $V^+$ to the set
\begin{equation}\notag
S^+=\{ (x,\eta) : x=(T_0, ... , T_{2N}), \ \rho\leq \eta\}.
\end{equation}

Vector fields $V^-$ and $V^+$ differ only in the expressions $\ov{a}_{2n}(\eta)$ \eqref{a2nbars} and $\ov{b}_{2n}(\eta)$ \eqref{b2nbars}. Note that $V^-$ is defined at points $(x,\rho)$, and further that
\begin{equation}\notag
\ov{a}_{2n}(\rho)=\al_2s_{2n}-(4n+1)(\al_2-\al_1)\int^\rho_0q_{2n}(y)dy)=\ov{b}_{2n}(\rho),
\end{equation}
that is, $V^-$ and $V^+$ agree on the set
\begin{equation}\label{Sigma}
\Sigma=\{ (x, \eta): x=(T_0, ... , T_{2N}), \ \eta=\rho\}.
\end{equation}
The set $\Sigma$ is known as a {\em switching boundary} (or {\em discontinuity boundary}) \cite{dibernardo} for the vector field
\begin{equation}\label{Vee}
V(x,\eta)=\begin{cases}
V^-(x,\eta),& (x,\eta)\in S^-\\
V^+(x,\eta),& (x,\eta)\in S^+.
\end{cases}
\end{equation}
The vector field $V$ is smooth on $\R^{N+2}\backslash \Sigma$ and continuous on $\R^{N+2}$. Standard ODE theory then implies solutions to the initial value problem
\begin{equation}\label{IVP}
(\dot{x},\dot{\eta})=V(x,\eta), \ (x(0),\eta(0))=(x_0,\eta_0)
\end{equation}
exist. However, $V$ fails to be $C^1$ on $\Sigma,$ as we now show.

Let $J^-$ and $J^+$ denote the Jacobian matrices of vector fields $V^-$ and $V^+$, respectively. 

\vspace{0.1in}
\noindent
{\bf Proposition 2.} \
For $(x,\rho)\in\Sigma, \ J^-(x,\rho)\neq J^+(x,\rho).$
\

\vspace{0.1in}
\noindent
{\em Proof}. Let $(x,\rho)\in \Sigma.$ Note
\begin{equation}\notag
\frac{\partial f^-_0}{\partial \eta}=\frac{Q}{B}(\al_i-\al_1)s(\eta) \ \mbox{ and } \ \frac{\partial f^+_0}{\partial \eta}=\frac{Q}{B}(\al_2-\al_1)s(\eta),
\end{equation}
expressions that do not agree when $\eta=\rho$, in particular. Hence $J^-(x,\rho)\neq J^+(x,\rho).$

\hfill \ $_\square$

\vspace{0.1in}
While $V$ is not $C^1$ at points in $\Sigma$, we do have that 
$V$ is locally Lipschitz, thereby guaranteeing that the model ODE \eqref{IVP} has unique  solutions.

\vspace{0.05in}
\noindent
{\bf Proposition 3.} \ 
{\em The vector field $V(x,\eta)$ given by system \eqref{Vee} is locally Lipschitz, that is, for any $z=(x,\eta)\in \R^{N+2}$, there exists a constant $L(z)$ and $\delta>0$ such that  for all $z_1, z_2\in B_\delta(z),$ $$\|V(z_1)-V(z_2)\|\leq L(z)\|z_1-z_2\|.$$
}

\vspace{0.05in}
\noindent
{\em Proof}. Each of $V^-$ and $V^+$ is smooth on $\R^{N+2}$, implying that $V^-$ and $V^+$ are each  locally Lipschitz on $\R^{N+2}$. It follows that $V$ is locally Lipschitz on $\R^{N+2}\backslash \Sigma.$ 

Recall that for all $(x,\rho)\in \Sigma, \ V^+(x,\rho)-V^-(x,\rho)=0$. By Proposition 2, $J^-(x,\rho)\neq J^+(x,\rho).$ In this scenario the switching  boundary $\Sigma$ is said to be {\em uniform with degree 2} (section 2.2 in \cite{dibernardo}). In particular, given $z_0=(x_0,\rho)\in\Sigma, $ there exists a $\delta^\prime>0$ and a smooth function $F:B_{\delta^\prime}(z_0)\rightarrow \R^{N+2}$ satisfying, for all $z=(x,\eta)\in  B_{\delta^\prime}(z_0),$
\begin{equation}\label{eff}
V^+(z)-V^-(z)=F(z)(\eta-\rho).
\end{equation}
This follows from consideration of the Taylor series  of the function $V^+-V^-$ expanded about $z_0=(x_0,\rho)$. (Note if we disregard connections to the model, each of $V^-$ and $V^+$ is defined on $\R^{N+2}$, and so each is defined on both sides of $\Sigma$.) The quantity   $\eta-\rho$ may be chosen as one of the coordinates in regions close to  $\Sigma$ and, since $V^+-V^-$ vanishes on $\Sigma$, 
 it  follows  that each term in the Taylor series expansion of each component function of $V^+-V^-$ has a factor of $\eta-\rho$. The expression $\eta-\rho$ is raised to the first power in equation \eqref{eff} since the Jacobian matrix of $V^+-V^-$ does not vanish on $\Sigma$ by Proposition 2.

 As each of $V^-$ and $V^+$ is locally Lipschitz at $z_0$, pick positive numbers $\delta_1, \delta_2, L_1$ and $L_2$ such that for all $z_1, z_2\in B_{\delta_1}(z_0)$,
\begin{equation}\notag
\|V^-(z_1)-V^-(z_2)\|\leq L_1\|z_1-z_2\|,
\end{equation}
and for all $z_1, z_2\in B_{\delta_2}(z_0)$,
\begin{equation}\notag
\|V^+(z_1)-V^+(z_2)\|\leq L_2\|z_1-z_2\|.
\end{equation}
Let $0<\delta<\min\{\delta^\prime, \delta_1, \delta_2\}$. Pick $K$ so that for all $z\in \ov{B_\delta(z_0)}, \ \|F(z)\|\leq K.$ Set $M=\max\{L_1, L_2\}+K$, and let $z_1, z_2\in B_\delta(z_0)$. If  $z_1, z_2\in S^-$ or $z_1, z_2\in S^+$, then clearly $\|V(z_1)-V(z_2)\|\leq M\|z_1-z_2\|.$

Without loss of generality, assume $z_1=(x_1,\eta_1)\in S^-$ and $z_2=(x_2,\eta_2)\in S^+$. Then
\begin{align}\notag
\|V(z_1)-V(z_2)\|&=\|V^-(z_1)-V^+(z_2)\|=\|V^-(z_1)-(V^-(z_2)+F(z_2)(\eta_2-\rho))\|\\ \notag
&\leq L_1\|z_1-z_2\|+K(\eta_2-\rho)\leq L_1\|z_1-z_2\|+K(\eta_2-\eta_1)\\ \notag
&\leq (L_1+K)\|z_1-z_2\|\leq M\|z_1-z_2\|,\notag
\end{align}
and we have that $V$ is locally Lipschitz on $\Sigma$. Combined with the fact $V$ is locally Lipschitz on $\R^{N+2}\backslash \Sigma$, we have the desired result. \hfill \ $_\square$

\subsection{Jormungand diffusion model dynamics}

We see via Proposition 3 that solutions to the nonsmooth system \eqref{Vee} exist and are unique. However, we cannot directly apply geometric singular perturbation theory to system \eqref{Vee} due to Proposition 2: Fenichel's Theorem applies to smooth vector fields (at least $C^1$), and system \eqref{Vee} is not smooth on the switching boundary $\Sigma$. As we are interested in a stable equilibrium solution with ice edge to the south of $y=\rho$, we first apply the theory  to system \eqref{Vminus}, subsequently restricting the domain of $\eta$ to return to the model interpretation. We repeat this process for $\eta>\rho$ (vector field \eqref{Vplus}), and then use Proposition 3 to analyze the behavior of system \eqref{Vee}. 

As in section 3, system \eqref{Vminus} has a globally attracting curve of rest points 
\begin{equation}\notag
\Lambda^0_-=\{(x,\eta) : x=h^0_-(\eta)\},  \ \ h^0_-(\eta)=(f^-_0(\eta), f^-_2(\eta), ... , f^-_{2N}(\eta)),
\end{equation}
when $\ep=0$. By Theorem 3.1, for $\ep$ sufficiently small system \eqref{Vminus} has an attracting invariant manifold on which the dynamics are well-approximated by the ODE
\begin{equation}\label{hminus}
\dot{\eta}=\ep\left(\sum^N_{n=0} f^-_{2n}(\eta)p_{2n}(\eta)-T_c\right)=\ep h^-(\eta).
\end{equation}
As each polynomial $f^-_{2n}$ has degree $2n+3$, we have that $h^-(\eta)$ is a polynomial of degree $4N+3$.

Set $\rho=0.35$ as in \cite{abb}. For the cold world of the Neoproterozoic ice ages, we let
\begin{equation}\label{Jparams}
Q=321, A=167, B=1.9, \al_1=0.32, \al_i=0.36, \al_2=0.8, T_c=0, D=0.25.
\end{equation}
 It is known that 600-700 million years ago the solar constant $Q$ was roughly 94\% of its current value of 343 W/m$^2$. Many other parameter values are taken from \cite{abb}, where the main criterion for the existence of a stable Jormungand state was found to be a significant difference between snow-covered ice albedo $\al_2$ and bare ice albedo $\al_i$. Relative to that used in section 3, we have chosen the smaller value $D=0.25$  as Abbot et al argue that meridional heat transport was less efficient during these periods of extensive glacial cover. Finally, the drawdown of atmospheric CO$_2$ via silicate weathering would have been greatly reduced due to the lack of significant ice free continental surfaces. As atmospheric CO$_2$ absorbs outgoing longwave radiation, we have reduced the $A+BT$-term in equation \eqref{budwiddiff}  by decreasing $A$.

\vspace{0.1in}
\noindent
{\bf Proposition 4.} \ 
{\em With parameters as in \eqref{Jparams} and $N=1$, system \eqref{Vee} has a stable Jormungand \eq solution for $\ep$ sufficiently small. That is, there is a stable \eq solution that does not correspond to a snowball Earth event and for which the ice line  lies within 20$^\circ$ of the equator. 
}

\vspace{0.1in}
\noindent
{\em Proof}. For sufficiently small $\ep$, the behavior of system \eqref{Vminus} is well-approximated by the dynamics of the ODE $\dot{\eta}=\ep h^-(\eta)$, where 
\begin{equation}\notag
 h^-(\eta)=f^-_0(\eta)p_0(\eta)+f^-_{2}(\eta)p_2(\eta),
 \end{equation}
by Theorem 3.1. One can   show the degree 7 polynomial $ h^-(\eta)$ has one real zero $\eta=\eta_1$ in $(0,\rho)$, and that $(h^-)^\prime(\eta_1)<0$ (see Figure 4(a)).  System \eqref{Vee} then has a stable equilibrium solution well-approximated by $(h^0_-(\eta_1),\eta_1)$; we note $y=\rho$ corresponds to 20.5$^\circ$N. \hfill $_\square$

\vspace*{.4in} 
\begin{figure}[t!]
\begin{center} \hspace*{0.325in}
\includegraphics[width=6in,trim = 1.6in 7.85in .6in  1in, clip]{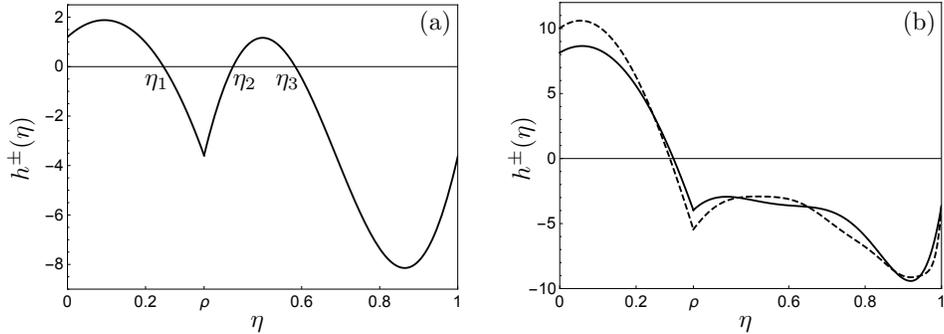}\\
\parbox{5in}{\caption{Plots of the Jormungand diffusion model functions $h^-(\eta) \, $   ($\eta<\rho) \, $ and $h^+(\eta) \, $   $ \, (\eta\geq \rho$) for $D=0.25$.
 \ (a) $ N=1$ in \eqref{TLeg}. \  (b)  {\em Solid}: $N=2$ in \eqref{TLeg}. {\em Dashed}: $N=5$ in \eqref{TLeg}.  }}
 
 \vspace*{.1in}
 \end{center}
\end{figure}

\vspace{-0.1in}
 One can apply Theorem 3.1 to system \eqref{Vplus} as well, and then restrict to the set $S^+$. For $\ep$ sufficiently small, system \eqref{Vplus} is  well-approximated by   the ODE $\dot{\eta}=\ep h^+(\eta)$, where 
\begin{equation}\notag
 h^+(\eta)=f^+_0(\eta)p_0(\eta)+f^+_{2}(\eta)p_2(\eta),
 \end{equation}
when $\rho<\eta$. We plot $h^+(\eta)$ in Figure 4 as well; for $N=1$ in equation \eqref{TLeg}, system \eqref{Vplus} has an unstable \eq solution with ice line near $\eta=\eta_2$ (intermediate ice cap), and a stable \eq solution with ice line near $\eta=\eta_3$ (small ice cap). This result is qualitatively the same  as that realized when using a smooth version of the Jormungand albedo function \eqref{Jalb} and Hadamard graph transform techniques (see Figure 8 in \cite{walshwid}). 

 Note the graph of the function 
\begin{equation}\notag
h(\eta)=
\begin{cases} h^-(\eta), & \eta<\rho\\ h^+(\eta), & \rho\leq\eta
\end{cases}
\end{equation}
is continuous at $\eta=\rho$.  As with the agreement of $V^-$ and $V^+$ on $\Sigma$, this follows from the fact $\ov{a}_{2n}(\eta)$ and $\ov{b}_{2n}(\eta)$ agree when $\eta=\rho$. The graphs of the invariant manifolds for vector fields \eqref{Vminus} and \eqref{Vplus}, however, are only within $O(\ep)$ of the graphs of $h^-$ and $h^+$ plotted in Figure 4(a). Nonetheless, system \eqref{Vee} does admit unique solutions by Proposition 3. Hence for $\ep$ sufficiently small, a trajectory of system \eqref{Vee} with $\eta(0)$ above $\eta=\rho$, but below the unstable equilibrium ice line $\eta=\eta_2$,  will see $\eta$ decrease as the trajectory passes through the switching boundary $\Sigma$. As seen in Figure 4(a), $h^-(\rho)<0$, and the $\eta$-value along this trajectory will continue to decrease, limiting on the large, stable ice line position near $\eta=\eta_1$.

Interestingly, only the large stable ice cap remains when higher-order approximations are used ($N=2, 5$ in equation \eqref{TLeg}; see Figure 4(b)). The incorporation of smaller spatial scale   information into the model makes for a colder world; in this case every trajectory converges to the stable Jormungand state.

\subsection{Bifurcations}

We note that several bifurcations occur in the Jormungand model with diffusive heat transport. With several parameters to choose from, the bifurcation parameter used here is the constant $A$ (recall the $A+BT$-term in equation \eqref{bud} models the outgoing longwave radiation). Atmospheric greenhouse gases such as CO$_2$ absorb OLR, so a decrease in $A$ might then be viewed as an increase in atmospheric CO$_2$, while an increase in $A$ can be associated with a decrease in CO$_2$. In this way the parameter $A$ can be viewed as a proxy for atmospheric greenhouse gas concentrations, although the rise and fall of the former corresponds to the fall and rise of the latter. 

In Figure 5 we plot the position of the ice line at equilibrium as the parameter $A$ is varied (with $N=5$ in \eqref{TLeg}). Not surprisingly, if $A$ is large enough (CO$_2$ concentrations are sufficiently low) the climate system heads to a snowball Earth as $\eta\rightarrow 0$. If $A$ is sufficiently small, so that OLR is absorbed to a great extent, the model tends toward an ice free Earth with $\eta\rightarrow 1$.

A first bifurcation occurs as $A$ increases through roughly 150, with the creation of an unstable \eq solution with a small ice cap, and a stable \eq solution with a slightly larger ice cap. This unstable solution is a manifestation of the SICI phenomenon mentioned above. 
For $A$-values less than roughly 157 there is no stable Jormungand state as every solution with $\eta(0)$ near 0 converges to the small stable ice line position.

For $A$ between roughly 157 and 161.5 the system exhibits bistability, with  a stable \eq solution with a moderately sized ice cap and a stable Jormungand state. As $A$ increases through 161.5 there is a saddle node bifurcation, leaving   the stable Jormungand state as the sole equilibrium solution. 

As $A$ continues to increase, so that the planet is radiating to space ever more efficiently, the stable ice line position heads to tropical latitudes. A final bifurcation occurs at $A\approx 187$, where the Jormungand state disappears in another saddle node bifurcation, with the  system now experiencing a runaway snowball Earth event.

It is interesting to note the bifurcation diagram plotted in Figure 5 exhibits  many of the characteristics of the bifurcation diagram for the full infinite-dimensional system \eqref{budwid} when using a smooth variant of the Jormungand albedo function \eqref{Jalb} (see Figure 9 in \cite{walshwid}). One exception is the appearance of  SICI in the Jormungand diffusion model; by contrast, in \cite{walshwid} the ice free Earth state is unstable in the sense  there is no equilibrium between the small stable ice line position and the North Pole. 

\begin{figure}[t!]
\begin{center}\hspace*{-.1in}
\includegraphics[width=4.8in,trim = 1.2in 6.7in 1.2in  1in, clip]{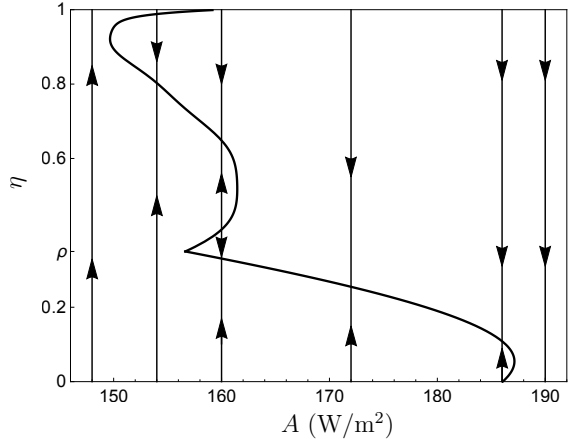}\\
\parbox{5in}{\caption{Bifurcation plot for the Jormungand model with diffusive heat transport, with $N=5$ in \eqref{TLeg}.}}
 \end{center}
\end{figure}

\section{Relaxation to the mean heat transport}

In this section we revisit and expand upon approximations of system \eqref{budwid} presented in \cite{dickesther} and \cite{walshrack}, investigating the effect (or  lack thereof) the inclusion of higher-order terms in equation \eqref{TLeg} has on results previously obtained.  We will also allow for the use of higher-order $s_{2n}$-terms in equation \eqref{fullsLeg}.

\subsection{Budyko's albedo function}

Consider system \eqref{budwid} with albedo function \eqref{budalb}; recall $\ov{T}$ denotes the global mean surface temperature. It can be shown that  equilibrium temperature profiles for uncoupled equation \eqref{budwidA} have the form
\begin{equation}\label{Tstarbud}
T^*(y,\eta)=\frac{1}{B+C}\left(Qs(y)(1-\al(y,\eta))-A+\frac{C}{B}(Q(1-\ov{\al})-A)\right),
\end{equation}
where $\ov{\al}=\int^1_0 \al(y)s(y) dy=\al_1 \int^\eta_0 s(y) dy + \al_2 \int^1_\eta s(y) dy$ (see \cite{dickesther}, e.g.). Assuming $s(y)$ as in \eqref{fullsLeg} is truncated at $n=N$, $T^*(y,\eta)$ is a piecewise degree $2N$ polynomial in $y$ with discontinuity at $y=\eta$, for each fixed $\eta$. Hence one might consider the use of Legendre polynomial approximations, given the form of the equilibrium solutions of equation \eqref{budwidA}, even though their role as eigenfunctions is no longer relevant. This idea was originally suggested in a preprint version of \cite{dickesther} in the case $N=1$, that is, when using quadratic approximations. Here we carry out the analysis for arbitrary $N$, while omitting   the more tedious computational details.

Express $T(y,t)$ piecewise as
\begin{equation}\label{TandV}
T(y,t)=\begin{cases}
U(y,t)=\sum^N_{n=0} T_{2n}(t)p_{2n}(y), & y<\eta\\
V(y,t)=\sum^N_{n=0} V_{2n}(t)p_{2n}(y), & \eta<y\\
\frac{1}{2}(U(\eta,t)+V(\eta,t)),& y=\eta.
\end{cases}
\end{equation}
The choice for $T(\eta,t)$ arises by taking the average of the limits of $U(y,t)$ and $V(y,t)$ as $y$ approaches $\eta$ from the left and right, respectively. 
Substituting $U(y,t)$ and $V(y,t)$ separately into equation \eqref{budwidA} yields
\begin{equation}\label{subUrelax}
R\sum^N_{n=0} \dot{T}_{2n}p_{2n}(y)=\sum^N_{n=0} ( Q(1-\al_1)s_{2n}-(B+C)T_{2n})p_{2n}(y)-A+C\ov{T}
\end{equation}
and
\begin{equation}\label{subVrelax}
R\sum^N_{n=0} \dot{V}_{2n}p_{2n}(y)=\sum^N_{n=0} ( Q(1-\al_2)s_{2n}-(B+C)V_{2n})p_{2n}(y)-A+C\ov{T},
\end{equation} 
where 
\begin{align}\label{relaxTbar1}
\ov{T}&=\int^\eta_0 U(y,t) dy+\int^1_\eta V(y,t)dy=\int^1_0 V(y,t)dy-\int^\eta_0 (V(y,t)-U(y,t)) dy\\\notag
&=V_0-\sum^N_{n=0} (V_{2n}-U_{2n})P_{2n}(\eta), \ \mbox{ with } \ P_{2n}(\eta)=\int^\eta_0 p_{2n}(y)dy.\notag
\end{align}
Equating coefficients of $p_{2n}$ in \eqref{subUrelax} and in \eqref{subVrelax} yields the system of ODE
\begin{align}\label{relaxUV1}
R\dot{T}_0&=Q(1-\al_1)-(B+C)T_0-A+C\ov{T},\\\notag
R\dot{T}_{2n}&=Qs_{2n}(1-\al_1)-(B+C)T_{2n}, \ n=1, ... , N,\\\notag
R\dot{V}_0&=Q(1-\al_2)-(B+C)V_0-A+C\ov{T},\\\notag
R\dot{V}_{2n}&=Qs_{2n}(1-\al_2)-(B+C)V_{2n}, \ n=1, ... , N.\notag
\end{align}
With an eye toward coupling system \eqref{relaxUV1} with the ice line equation \eqref{budwidB}, we note
\begin{equation}\label{Tofeta1}
T(\eta,t)=\frac{1}{2}\sum^N_{n=0}(T_{2n}+V_{2n})p_{2n}(y).
\end{equation} 
The change of variables $u=\frac{1}{2}(T_0+V_0), \ v=T_0-V_0$ proves helpful, resulting in the system
\begin{subequations}\label{relaxUV2}
\begin{eqnarray}
R\dot{u}& \hspace*{-.075in} =&  \hspace*{-.075in} Q(1-\al_0)-(B+C)u-A+C\ov{T},\label{relaxUV2A}\\
R\dot{T}_{2n} & \hspace*{-.075in} = &  \hspace*{-.075in} Qs_{2n}(1-\al_1)-(B+C)T_{2n}, \ n=1, ... , N,\label{relaxUV2B}\\
R\dot{v} & \hspace*{-.075in} = &  \hspace*{-.075in} Q(\al_2-\al_1)-(B+C)v,\label{relaxUV2C}\\
R\dot{V}_{2n} & \hspace*{-.075in} = &  \hspace*{-.075in} Qs_{2n}(1-\al_2)-(B+C)V_{2n}, \ n=1, ... , N,\label{relaxUV2D}
\end{eqnarray}\end{subequations}
where $\al_0=\frac{1}{2}(\al_1+\al_2)$. In the new variables,
\begin{equation}\label{relaxTbar2}
\ov{T}=
u+v(\eta-{\textstyle{\frac{1}{2}}})-\sum^N_{n=1} (V_{2n}-U_{2n})P_{2n}(\eta), \notag
\end{equation}
and
\begin{equation}\label{Tofeta2}
T(\eta,t)=u+\frac{1}{2}\sum^N_{n=1}(T_{2n}+V_{2n})p_{2n}(y).
\end{equation}
Set $L=Q/(B+C).$ Note that as $t\rightarrow\infty$,
\begin{equation}\label{equils}
T_{2n}\rightarrow T^*_{2n}=Ls_{2n}(1-\al_1), \ V_{2n}\rightarrow V^*_{2n}=Ls_{2n}(1-\al_2), \ n=1, ... , N, \end{equation}
\begin{equation}\notag 
\text{ and } \ 
v\rightarrow v^*=L(\al_2-\al_1).
\end{equation}
 Hence the long term behavior of system \eqref{relaxUV2} reduces to that of equation  \eqref{relaxUV2A} when  assuming variables other than $u$ and $\eta$ are at equilibrium. Coupled with equation \eqref{budwidB}, the resulting two equations can be placed in the form
\begin{align}\label{ueta}
\dot{u}&=-\frac{B}{R}(u-f(\eta)),\\\notag
\dot{\eta}&=\ep\left(u+\frac{1}{2}\sum^N_{n=1}(T^*_{2n}+V^*_{2n})p_{2n}(\eta)-T_c\right),\notag
\end{align}
where
\begin{align}\notag
f(\eta)&=\frac{1}{B}\left(Q(1-\al_0)-A+C\left(v^*(\eta-{\textstyle{\frac{1}{2}}})+\sum^N_{n=1}(V^*_{2n}-T^*_{2n})P_{2n}(\eta)\right)\right)\\\notag
&=\frac{1}{B}\left(Q(1-\al_0)-A+CL(\al_2-\al_1)\left(\eta-{\textstyle{\frac{1}{2}}}-\sum^N_{n=1}s_{2n}P_{2n}(\eta)\right)\right).
\end{align}
When $\ep=0$ system \eqref{ueta} has a globally attracting curve of rest points $(u,\eta)=(f(\eta),\eta)$. By Theorem 3.1, system \eqref{ueta} is well-approximated by the ODE
\begin{align}\label{hrelaxbud}
\dot{\eta}&=\ep\left(f(\eta)+\frac{1}{2}\sum^N_{n=1}(T^*_{2n}+V^*_{2n})p_{2n}(\eta)-T_c\right)\\ \notag
&=\ep\left(f(\eta)+L(1-\al_0)\sum^N_{n=1}s_{2n}p_{2n}(\eta) \, -T_c\right)\\ \notag
&=\ep h(\eta),
\end{align}
for $\ep$ sufficiently small. In general, $h$ is a polynomial in $\eta$ of degree $2N+1$.

In the case $N=1$ the model development and analysis above reduces to that presented in section 3.2 in \cite{walshrack}. Using quadratic approximations for $s(y)$ and piecewise for $T(y,t)$ as well, the function $h(\eta)$ given by equation \eqref{hrelaxbud} is a cubic polynomial with two zeros in [0,1] when using  parameter values \eqref{budparams} (and $C=3.09$). One is a stable equilibrium with $\eta$ near $\eta_2\approx 0.94$,  while the other is  an unstable \eq with $\eta$ near $\eta_1\approx 0.24$ (see Figure 6(a)).

\begin{figure}[b!]
\begin{center} \hspace*{0.2in}
\includegraphics[width=6in,trim = 1.7in 7.8in .65in  1in, clip]{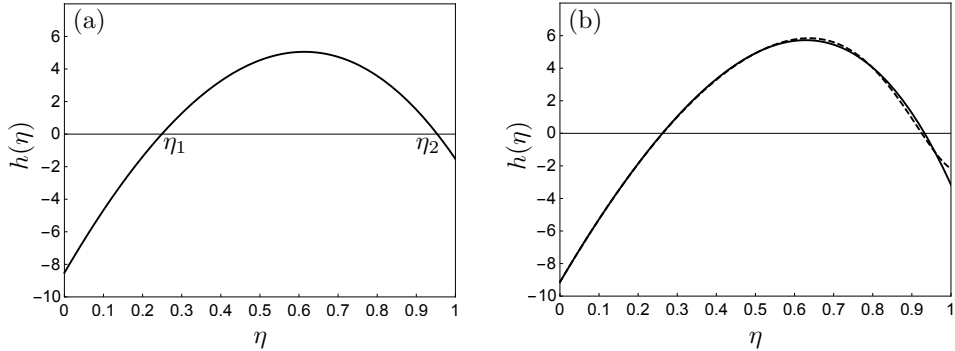}\\
\parbox{5in}{\caption{Plots of function $h(\eta)$ given by \eqref{hrelaxbud}. \ (a) $ N=1$ in \eqref{TLeg}. \  (b) Including higher-order terms in \eqref{fullsLeg}. {\em Solid}: $N=2$  ($s_4=-0.044$). {\em Dashed}: $N=5$ ($s_6=0.006, s_8=0.016, s_{10}=0.006$). Parameters as in \eqref{budparams}, $C=3.09$.}}
 \end{center}
\end{figure}

If we restrict to the use of quadratic polynomials in the approximations, the behaviors of the diffusion and relaxation to the mean models are qualitatively the same for appropriately chosen $D$-values (compare Figures 6(a) and 2(b)). This is not surprising at first glance, for if $\eta$ is kept fixed and equation \eqref{TLeg} with $N=1$ is substituted into equation \eqref{diffuse}, the result is  $-6DT_2(t)p_2(y)$. Note in the $N=1$ case $T(y,t)=T_0(t)+T_2(t)p_2(y)$, so that 
\begin{equation}\notag
6DT_2(t)p_2(y)=6D(T(y,t)-T_0(t))=6D(T(y,t)-\ov{T}).
\end{equation}
We see that for a second mode approximation, and with the ice line fixed, the heat transport in the diffusion model is treated the same as that in  relaxation to the mean model, varying only by the choice of constants $D$ and $C$ (see \cite{north}).

We further note that if we set $s_{2n}=0$ for $n>1$ as in previous sections, $T_{2n}\rightarrow 0$ and $V_{2n}\rightarrow 0$ as $n\rightarrow \infty$ for $n>1$, implying a return to the model behavior when using a quadratic approximation for $U(y,t)$ and $V(y,t)$ (and so for $T(y,t)$).

The only remaining question pertains to the use of higher order terms in expression \eqref{fullsLeg} for $s(y)$. The $s_{2n}$-terms appear only in the equilibrium point expressions \eqref{equils}. That is, $T_{2n}$ and $V_{2n}$ still converge to equilibrium points but, if $s_{2n}\neq0$ for $n>1$, $T^*_{2n}$ and $V^*_{2n}$ are no longer 0. This in turn has an effect on $h(\eta)$ in equation \eqref{hrelaxbud}. This effect is negligible, however, as can be seen in Figure 6(b). The qualitative analysis of the corresponding ODE \eqref{hrelaxbud} remains unchanged from the case in which $s_{2n}=0$ for $n>1$. The analysis presented in \cite{walshrack} using quadratic polynomials, following that presented in \cite{dickesther}, is sufficient to understand any approximation to  system 
\eqref{budwid} via the use of Legendre polynomials.

\subsection{The Jormungand albedo function}

In this section we consider approximations to system \eqref{budwid} in the case where the albedo function is given by \eqref{Jalb}.  We seek to establish the existence of a stable \eq solution with a large ice cap, as was the case with the diffusion model in section 4. Note when $\eta>\rho$, albedo functions \eqref{budalb} and \eqref{Jalb} are identical, and hence the analysis follows that presented in the previous section. However, we now let $h^+(\eta)$ denote the function previously defined as $h(\eta)$ in equation \eqref{hrelaxbud}, the plus sign indicating the restriction of $\eta>\rho$ in the Jormungand setting. The plot of $h^+(\eta)$ can be seen in Figure 7.

For $\eta<\rho$, express $T(y,t)$ piecewise as
\begin{equation}\label{TVW}
T(y,t)=\begin{cases}
U(y,t)=\sum^N_{n=0} T_{2n}(t)p_{2n}(y), & y<\eta\\
V(y,t)=\sum^N_{n=0} V_{2n}(t)p_{2n}(y), & \eta<y<\rho\\
W(y,t)=\sum^N_{n=0} W_{2n}(t)p_{2n}(y), & \rho<y.
\end{cases}
\end{equation}
Set $T(\eta,t)=\frac{1}{2}(U(\eta,t)+V(\eta,t)),$ and   $T(\rho,t)=\frac{1}{2}(V(\rho,t)+W(\rho,t)).$  

Details for the following derivation of the model equations can be found in the Appendix. Plugging $U, V$ and $W$ separately into equation \eqref{budwidA} and equating coefficients of $p_{2n}(y)$ yields
\begin{align}\label{relaxTVW1}
R\dot{T}_0&=Q(1-\al_1)-(B+C)T_0-A+C\ov{T},\\\notag
R\dot{T}_{2n}&=Qs_{2n}(1-\al_1)-(B+C)T_{2n}, \ n=1, ... , N,\\\notag
R\dot{V}_0&=Q(1-\al_i)-(B+C)V_0-A+C\ov{T},\\\notag
R\dot{V}_{2n}&=Qs_{2n}(1-\al_i)-(B+C)V_{2n}, \ n=1, ... , N.\\\notag
R\dot{W}_0&=Q(1-\al_2)-(B+C)W_0-A+C\ov{T},\\\notag
R\dot{W}_{2n}&=Qs_{2n}(1-\al_2)-(B+C)W_{2n}, \ n=1, ... , N.\notag
\end{align}
Recalling $P_{2n}(\eta)=\int^\eta_0 p_{2n}(y)dy$ and simplifying the expression for $\ov{T}$, one finds
\begin{align}\label{TbarTVW1}
\ov{T}&=\int^\eta_0U(y,t) dy+\int^\rho_\eta V(y,t)dy+\int^1_\rho W(y,t) dy\\\notag
&=\sum^N_{n=0}(T_{2n}-V_{2n})P_{2n}(\eta)+\sum^N_{n=0}(V_{2n}-W_{2n})P_{2n}(\rho).\notag
\end{align}
We also have 
\begin{equation}\label{JTofeta1}
T(\eta,t)=\frac{1}{2} \ \sum^N_{n=0}(T_{2n}+V_{2n})p_{2n}(\eta).
\end{equation}
A sequence of variable changes converts system \eqref{relaxTVW1} into
\begin{subequations}\label{TVWlast}
\begin{eqnarray}
R\dot{w}& \hspace*{-.075in} =&  \hspace*{-.075in} Q(1-\gamma_1)-A-(B+C)w+C\ov{T}, \label{TVWlastA}\\
R\dot{z}_2& \hspace*{-.075in} =&  \hspace*{-.075in} Q(\al_i-\al_0)-(B+C)z_2, \label{TVWlastB}\\
R\dot{z}_1& \hspace*{-.075in} =&  \hspace*{-.075in} Q(\al_2-\al_1)-(B+C)z_1, \label{TVWlastC}\\
R\dot{T}_{2n}& \hspace*{-.075in} =&  \hspace*{-.075in} Qs_{2n}(1-\al_1) -(B+C)T_{2n}, \ n=1,2, .. , N, \label{TVWlastD}\\
R\dot{V}_{2n}& \hspace*{-.075in} =&  \hspace*{-.075in} Qs_{2n}(1-\al_i) -(B+C)V_{2n}, \ n=1,2, .. , N, \label{TVWlastE} \\
R\dot{W}_{2n}& \hspace*{-.075in} =&  \hspace*{-.075in} Qs_{2n}(1-\al_2) -(B+C)W_{2n}, \ n=1,2, .. , N, \label{TVWlastF}
\end{eqnarray}
\end{subequations}
where $\al_0=\frac{1}{2}(\al_1+\al_2)$ and $\gamma_1=\frac{1}{2}(\al_0+\al_i)$. In the new variables,
\begin{align}\label{JTbarlast}
\ov{T}&=w+{\textstyle{\frac{1}{2}}}z_1(\eta+\rho-1)+z_2(\eta-\rho+{\textstyle{\frac{1}{2}}}) +\\\notag
& \hspace*{.25in} \sum^N_{n=1}(T_{2n}-V_{2n})P_{2n}(\eta)+\sum^N_{n=1}(V_{2n}-W_{2n})P_{2n}(\rho),\notag
\end{align}
and
\begin{equation}\label{JTofetalast}
T(\eta,t)=w+{\textstyle{\frac{1}{4}}}z_1+{\textstyle{\frac{1}{2}}}\sum^N_{n=1}(T_{2n}+V_{2n})p_{2n}(\eta).
\end{equation}
System \eqref{TVWlast} decouples nicely, and furthermore
\begin{equation}\label{Jequils}
z_2\rightarrow z^*_2=L(\al_i-\al_0), \ z_1\rightarrow z^*_1=L(\al_2-\al_1), \ T_{2n}\rightarrow T^*_{2n}=Ls_{2n}(1-\al_1),
\end{equation}
\begin{equation}\notag
V_{2n}\rightarrow V^*_{2n}=Ls_{2n}(1-\al_i),  \ \mbox{ and } \ W_{2n}\rightarrow W^*_{2n}=Ls_{2n}(1-\al_2)
\end{equation}
as $t\rightarrow\infty, \ n=1, ... , N$. Hence the behavior of system \eqref{TVWlast} reduces to that of equation \eqref{TVWlastA} when assuming all variables other than $w$ and $\eta$ are at equilibrium. Coupled with equation \eqref{budwidB}, and evaluating equations \eqref{JTbarlast} and \eqref{JTofetalast} at equilibrium values \eqref{Jequils}, the resulting  $(w,\eta)$-equations can be put into the form
\begin{align}\label{weta}
\dot{w}&=-\frac{B}{R}(w-f(\eta)),\\\notag
\dot{\eta}&=\ep\left(w+{\textstyle{\frac{1}{4}}}z^*_1+{\textstyle{\frac{1}{2}}}\sum^N_{n=1}(T^*_{2n}+V^*_{2n})p_{2n}(\eta)-T_c\right),\notag
\end{align}
where $f(\eta)=\frac{1}{B}(Q(1-\gamma_1)-A+Cg(\eta))$, with  $g(\eta)$  given by
\begin{align}\label{gofeta}
g(\eta)&={\textstyle{\frac{1}{2}}}z^*_1(\eta+\rho-1)+z^*_2(\eta-\rho+{\textstyle{\frac{1}{2}}}) +\\\notag
& \hspace*{.25in} \sum^N_{n=1}(T^*_{2n}-V^*_{2n})P_{2n}(\eta)+\sum^N_{n=1}(V^*_{2n}-W^*_{2n})P_{2n}(\rho).\notag
\end{align}
Substituting expressions \eqref{Jequils} into \eqref{weta} and simplifying, one arrives at
\begin{align}\label{newweta}
\dot{w}&=-\frac{B}{R}(w-f(\eta)),\\\notag
\dot{\eta}&=\ep\left(w+{\textstyle{\frac{1}{4}}}L(\al_2-\al_1)+L(1-\gamma_2)\sum^N_{n=1}s_{2n}p_{2n}(\eta) \, -T_c\right),\notag
\end{align}
where $\gamma_2=\frac{1}{2}(\al_1+\al_i), \, f(\eta) $ is as in \eqref{weta}, and \eqref{gofeta} becomes
\begin{align}\notag
g(\eta)&=-{\textstyle{\frac{1}{4}}}L(3\al_2-2\al_i-\al_1)+L(\al_i-\al_1)\sum^N_{n=0}s_{2n}P_{2n}(\eta)\\\notag
&\hspace*{.25in}+L(\al_2-\al_i)\sum^N_{n=0}s_{2n}P_{2n}(\rho).
\end{align}
We remark that $g(\eta)$, and hence $f(\eta)$, is a degree $2N+1$ polynomial in $\eta$.

There is a globally attracting curve of rest points $(w,\eta)=(f(\eta),\eta)$ when $\ep=0$ in \eqref{weta}. Invoking Fenichel's Theorem, the behavior of system \eqref{weta} can be discerned by consideration of the ODE
\begin{equation}\label{Jrelaxhminus}
\dot{\eta}=\ep\left(f(\eta)+{\textstyle{\frac{1}{4}}}L(\al_2-\al_1)+L(1-\gamma_2)\sum^N_{n=1}s_{2n}p_{2n}(\eta) \, -T_c\right)=\ep h^-(\eta),
\end{equation}
for $\ep$ sufficiently small. 

We have that system \eqref{relaxTVW1} is well-approximated by the ODE $\dot{\eta}=\ep h^-(\eta)$ when $\eta<\rho$ and by  the ODE $\dot{\eta}=\ep h^+(\eta)$ when $\eta>\rho$. Notably, $h^-(\rho)\neq h^+(\rho)$, so that the analysis leads to a discontinuous differential equation. The discontinuity at $\eta=\rho$ arises from the fact the temperature at the ice line, $T(\eta,t)$, is not continuous at $\eta=\rho$: for $\eta<\rho, T(\eta,t)=(U(\eta,t)+V(\eta,t))/2$, while $T(\eta,t)=(U(\eta,t)+W(\eta,t))/2$ when $\eta>\rho$.

We are lead to consider the one-dimensional differential inclusion
\begin{equation}\label{filippov}
\dot{\eta}\in h(\eta)= \begin{cases} h^-(\eta), & \eta<\rho\\
\{ (1-q)(h^-)^\prime(\rho)+q(h^+)^\prime(\rho) : q\in [0,1]\}, & \eta=\rho \\
h^+(\eta), & \rho<\eta,
\end{cases}
\end{equation}
to which the theory of Filippov \cite{fil} can be applied. The analysis of this Filippov flow in the case $N=1$ was presented in section 4.3 in \cite{walshrack}, to which we refer the reader for details. 

In Figures 7(a) and 7(b), parameters were chosen as in Figure 13 in \cite{walshrack}. Referring to Figures 7(a) and 7(b),   system \eqref{filippov} has a stable \eq solution with ice line near $\eta=\eta_3$, and an unstable \eq solution with ice line near $\eta=\eta_2$, for sufficiently small $\ep$. We note a similar $(\eta_2,\eta_3)$ unstable-stable pair arose in the $N=1$ case in the Jormungand diffusion model  (recall Figure 4(a)).

\begin{figure}[t!]
\begin{center}\hspace*{0.25in}
\includegraphics[width=6in,trim = 1.7in 6.54in .7in  1in, clip]{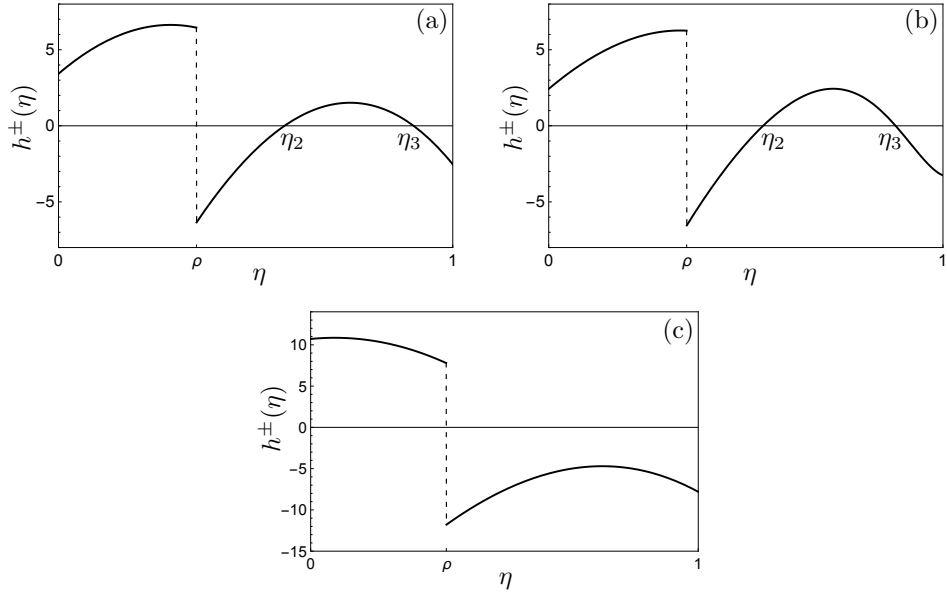}\\
\parbox{5in}{\caption{Plots of functions  $h^-(\eta)$ \eqref{Jrelaxhminus} (for $\eta<\rho$) and $h^+(\eta)$ (for $\eta> \rho$) 
(a) $ N=1$  in \eqref{TLeg}. \  (b) $ N=5$ ($s_4, ...  s_{10}$ as in Figure 6). Parameters for (a) and (b) as in Figure 13 in \cite{walshrack}. \  (c) $N=1$, parameters as in \eqref{Jparams}.   }}
 \end{center}
\end{figure}

Referring to the plot of $h^-(\eta)$ and $h^+(\eta)$ in Figure 7(a),
a solution with $\eta(0)\in (0,\rho)$ will see $\eta$ increase and reach  $\eta=\rho$ in finite time. A solution with $\eta(0)$ just to the right of $\eta=\rho$ will see $\eta$ decrease and again reach  $\eta=\rho$ in finite time. In this case the stable Jormungand state has the ice line sitting at $\eta=\rho$.

In Figure 7(c) we use the parameter values given in \eqref{Jparams}. In this case every solution has $\eta(t)$ reaching $\eta=\rho$ in finite time. This behavior is the discontinuous analog of the behavior of trajectories  in the diffusion model when higher-order terms in \eqref{TLeg}  were included in the analysis (recall Figure 4(b)).

If  we set $s_{2n}=0$ for $n>1$, then each of $T_{2n}, V_{2n}$ and $W_{2n}$ approaches 0 over time for $n=1, ... , N$. Hence the long term behavior of system \eqref{relaxTVW1} in this case is equivalent to the $N=1$ case (quadratic approximations). As  with Budyko's albedo function in section 5.1, including higher-order terms from expression \eqref{fullsLeg} in the analysis does not alter the model behavior in any qualitative sense (see Figure 7(b)). The use of quadratic polynomials is sufficient when approximating the temperature-ice line system \eqref{budwid} with Legendre polynomials when using the Jormungand albedo function.

\section{Conclusion}

The energy balance  models introduced by Budyko \cite{bud} and Sellers \cite{sell} were used to investigate the  role positive ice albedo feedback plays in climate dynamics. Budyko used a relaxation to the mean horizontal heat transport model component and the straightforward albedo function \eqref{budalb}. E. Widiasih enhanced the model by coupling Budyko's temperature equation to a dynamic ice line \cite{wid}. Using a smooth variant of albedo function \eqref{budalb} and working in the  function space setting, Widiasih proved the existence of a stable equilibrium with a small ice cap, and an unstable equilibrium with a large ice cap, provided the time constant $\ep$  for the evolution of the ice line met certain conditions. In particular, she found that  positive ice albedo feedback did not lead to an ice free Earth, although this feedback did lead to a snowball Earth if the ice line were ever positioned sufficiently far south. 

By adjusting the smooth albedo function, Walsh and Widiasih proved that a stable equilibrium solution with the ice line in tropical latitudes  could be realized in the coupled   model in the function space setting \cite{walshwid}, again assuming certain restrictions on $\ep$. This stable ``Jormungand" climate state had previously been realized using an idealized GCM in \cite{abb}, a work in which the aim  was to model certain extensive glacial episodes from Earth's past.

A second approach to studying the temperature-ice line model was presented in \cite{dickesther}, in the case where the albedo  was   given by the piecewise constant function \eqref{budalb} (and still incorporating  relaxation to the mean meridional heat transport). McGehee and Widiasih used quadratic polynomial approximations for the temperature profile and deduced  that the corresponding  system dynamics reduced to that of  a pair of ODE. Once again an unstable-stable pair of equilibria were found, similar in spirit to the result of Widiasih in the function space setting.

A similar quadratic approximation approach was used in \cite{walshrack}, albeit with piecewise constant albedo functions (including \eqref{Jalb}) modeled after an  albedo function used in \cite{abb}. Although tools from nonsmooth dynamical systems theory had to be appealed to in the analysis, the results were akin to those found in \cite{walshwid}, namely, that the system arrived at via quadratic approximations exhibited a stable Jormungand climate state for $\ep$ sufficiently small. 

In this paper we have replaced the relaxation to the mean heat transport with a diffusive process. Many authors have previously analyzed Budyko's equation with diffusive heat transport, in each case, however,  with the ice line fixed.  Due to the form of the spherical diffusion operator, approximations with Legendre polynomials is a natural tool to use in the analysis. We thus coupled the temperature equation with the ice line equation and applied the spectral method,    with diffusion constant $D$.
As  we assume the ice line moves slowly relative to changes in the temperature, the main instrument used to analyze the approximating systems of ODE is geometric singular perturbation theory.

When using Budyko's albedo function \eqref{budalb} we find the diffusive model, in essence,  recreates the dynamics seen in \cite{wid} and \cite{dickesther} for a range of $D$-values. If the diffusive transport is too efficient ($D$ is too high), however, the small stable ice cap is lost and a retreating ice line tends to the North Pole. For other ranges of $D$-values, a small stable \eq ice cap and an even smaller unstable \eq ice cap coexist. This small ice cap instability (SICI), previously known to occur in Budyko's (uncoupled) model, can therefore arise in the coupled model as well.  Incorporating  modes of order greater than 2 in the analysis
did not change the qualitative nature of model behavior.

We have applied the spectral method to the coupled model with diffusive transport using the  piecewise constant Jormungand albedo function \eqref{Jalb}. While the resulting vector field is nonsmooth due to the existence of a switching boundary, the vector field has been shown to be  locally Lipschitz.   This spectral approach produces a stable Jormungand state for a range of $D$-values, akin to the results in \cite{walshwid} and \cite{walshrack}. Incorporating higher modes into the approximation, however, changes the dynamics in the Jormungand model. For example, the model can exhibit  bistability when using a 2-mode approximation, while the stable Jormungand state becomes the sole equilibrium solution when using a 4-mode approximation.

Each of the results for the coupled model with diffusive heat transport were arrived at via geometric singular perturbation theory and, as such, $\ep$ is assumed to be sufficiently small.

For completeness we revisited both the Budyko and Jormungand coupled models in the case of relaxation to the mean heat transport, albeit using higher mode approximations than those used in \cite{dickesther} and \cite{walshrack}. In every case we found approximations by quadratic polynomials was sufficient to discern the  behavior of solutions to the corresponding systems of ODE.

We view this work  as completing the analysis of the dynamical system that, originating with   Budyko's seminal surface temperature equation \cite{bud} and its coupling with Widiasih's dynamic ice line equation \cite{wid},  was further analyzed  in \cite{walshwid}, \cite{dickesther} and \cite{walshrack}.

\vspace*{0.2in}
\noindent
{\bf Acknowledgment.} The author recognizes and appreciates the support of the Mathematics and Climate Research Network (http://www.mathclimate.org).

\vspace*{0.2in}
\noindent
{\bf Appendix.} \ We derive model equations \eqref{TVWlast}. Substituting $U, V$and $W$ from \eqref{TVW} into \eqref{budwidA} yields the three equations
\begin{align}\notag
R\sum^N_{n=0} \dot{T}_{2n}p_{2n}(y)&=\sum^N_{n=0} ( Q(1-\al_1)s_{2n}-(B+C)T_{2n})p_{2n}(y)-A+C\ov{T},\\ \notag
R\sum^N_{n=0} \dot{V}_{2n}p_{2n}(y)&=\sum^N_{n=0} ( Q(1-\al_i)s_{2n}-(B+C)V_{2n})p_{2n}(y)-A+C\ov{T},\\ \notag
R\sum^N_{n=0} \dot{W}_{2n}p_{2n}(y)&=\sum^N_{n=0} ( Q(1-\al_2)s_{2n}-(B+C)W_{2n})p_{2n}(y)-A+C\ov{T}.\notag
\end{align}
Equating coefficients of $p_{2n}(y)$ in each of the above equations then gives system \eqref{relaxTVW1}. The associated expression $\ov{T}$ in \eqref{TbarTVW1} is a function of $T_{2n}, V_{2n}, W_{2n}$ and $\eta$, while $T(\eta,t)$ is given by \eqref{JTofeta1}.

With an eye toward reducing  the number of occurrences of $\ov{T}$  in \eqref{relaxTVW1}, consider the change of variables 
\begin{equation}\notag
x={\textstyle{\frac{1}{2}}}(T_0+W_0), \ z_1=T_0-W_0.
\end{equation}
System \eqref{relaxTVW1} becomes
\begin{align}\label{TVWapp2}
R\dot{x}&= Q(1-\al_0)-(B+C)x-A+C\ov{T}  \\  \notag
R\dot{z}_1&= Q(\al_2-\al_1)-(B+C)z_1  \\ \notag
R\dot{T}_{2n}&  = Qs_2(1-\al_1) -(B+C)T_{2n}, \, n=1, ... , N,  \\  \notag
R\dot{V}_0&= Q(1-\al_i)-(B+C)V_0-A+C\ov{T}  \\  \notag
R\dot{V}_{2n}&  = Qs_2(1-\al_i) -(B+C)V_{2n}, \, n=1, ... , N, \\  \notag
R\dot{W}_{2n}&= Qs_2(1-\al_2) -(B+C)W_{2n}, \, n=1, ... , N, \notag
\end{align}
where $\al_0=\frac{1}{2}(\al_w+\al_s)$. Substituting $T_0=x+\frac{1}{2}z_1$ and $W_0=x-\frac{1}{2}z_1$ into \eqref{TbarTVW1} and \eqref{JTofeta1}, one finds
\begin{align}\label{Tbarapp2}
\ov{T}&= (x+{\textstyle{\frac{1}{2}}}z_1)\eta+(x-{\textstyle{\frac{1}{2}}}z_1)(1-\rho)+V_0(\rho-\eta) +\sum^N_{n=1}(T_{2n}-V_{2n})P_{2n}(\eta)\\\notag
& \hspace*{.5in}+ \sum^N_{n=1}(V_{2n}-W_{2n})P_{2n}(\rho),\notag
\end{align}
and 
\begin{equation}\label{Tofetaapp2}
T(\eta,t)={\textstyle{\frac{1}{2}}}(x+{\textstyle{\frac{1}{2}}}z_1+V_0)+{\textstyle{\frac{1}{2}}}\sum^N_{n=1}(T_{2n}+V_{2n})p_{2n}(\eta).
\end{equation}
For the second change of variables, let
\begin{equation}\notag
w={\textstyle{\frac{1}{2}}}(x+V_0), \ z_2=x-V_0.
\end{equation}
System \eqref{TVWapp2} becomes
\begin{align}\label{TVWapp3}
R\dot{w}&=Q(1-\gamma_1)-(B+C)w-A+C\ov{T},\\ \notag
R\dot{z}_2&=Q(\al_i-\al_0)-(B+C)z_2,\\ \notag
R\dot{z}_1&= Q(\al_2-\al_1)-(B+C)z_1  \\ \notag
R\dot{T}_{2n}&  = Qs_2(1-\al_1) -(B+C)T_{2n}, \, n=1, ... , N,  \\  \notag
R\dot{V}_{2n}&  = Qs_2(1-\al_i) -(B+C)V_{2n}, \, n=1, ... , N, \\  \notag
R\dot{W}_{2n}&= Qs_2(1-\al_2) -(B+C)W_{2n}, \, n=1, ... , N, \notag
\end{align}
where $\gamma_1=\frac{1}{2}(\al_0+\al_i)$. Substituting $x=w+\frac{1}{2}z_2$ and $V_0=w-\frac{1}{2}z_2$ into \eqref{Tbarapp2} and \eqref{Tofetaapp2}, one finds
\begin{align}\notag
\ov{T}&= w+ {\textstyle{\frac{1}{2}}}z_1(\eta+\rho-1)+z_2(\eta-\rho+{\textstyle{\frac{1}{2}}})+
\sum^N_{n=1}(T_{2n}-V_{2n})P_{2n}(\eta)\\\notag
& \hspace*{.5in}+ \sum^N_{n=1}(V_{2n}-W_{2n})P_{2n}(\rho)\notag
\end{align}
as in \eqref{JTbarlast}, and 
\begin{equation}\notag
T(\eta,t)=w+{\textstyle{\frac{1}{4}}}z_1+{\textstyle{\frac{1}{2}}}\sum^N_{n=1}(T_{2n}+V_{2n})p_{2n}(\eta),
\end{equation}
as in \eqref{JTofetalast}.

\vspace*{0.2in}

\renewcommand{\refname}

\end{document}